\documentclass[a4paper,10pt]{article}
\usepackage{algorithm2e}
\usepackage[latin1]{inputenc}
\usepackage[T1]{fontenc}
\usepackage{bezier}
\usepackage{dsfont}
\usepackage{mathptmx}      
\usepackage{amsfonts}
\usepackage{amsmath}
\usepackage{amssymb}
\usepackage{amsbsy}
\usepackage{amscd}
\usepackage{amsopn}
\usepackage{amssymb}
\usepackage{amstext}
\usepackage{amsxtra}
\usepackage{graphicx,psfrag}
\usepackage{latexsym}
\usepackage{multirow}
\newtheorem{thm}{Theorem}
\newtheorem{lem}{Lemma}

\def\squareforqed{\hbox{\rlap{$\sqcap$}$\sqcup$}}
\def\qed{\ifmmode\else\unskip\quad\fi\squareforqed}
\def\smartqed{\def\qed{\ifmmode\squareforqed\else{\unskip\nobreak\hfil
\penalty50\hskip1em\null\nobreak\hfil\squareforqed
\parfillskip=0pt\finalhyphendemerits=0\endgraf}\fi}}
\smartqed
\def\trigforqed{\hbox{\rlap{$\blacktriangle$}}}
\def\qend{\ifmmode\else\unskip\quad\fi\trigforqed}
\def\smartqned{\def\qend{\ifmmode\trigforqed\else{\unskip\nobreak\hfil
\penalty50\hskip1em\null\nobreak\hfil\trigforqed
\parfillskip=0pt\finalhyphendemerits=0\endgraf}\fi}}
\smartqned

\DeclareMathAlphabet{\mathcal}{OMS}{cmsy}{m}{n}

\newcommand{\wend}[1]{\smartqned}
\newcommand{\wabs}[1]{\left|#1\right|}

\newcommand{\wcal}[1]{\mathcal{#1}}

\newcommand{\wdsf}[2]{{#1}''\!\!\left(#2\right)}

\newcommand{\werror}[1]{\wcal{E}_{#1}}

\newcommand{\wf}[2]{{#1}\left(#2\right)}
\newcommand{\wfc}[2]{{#1}\!\left(#2\right)}
\newcommand{\wfcc}[2]{{#1}\!\!\left(#2\right)}

\newcommand{\wfloor}[1]{\lfloor {{#1}} \rfloor }

\newcommand{\wi}[1]{\mathds{I}}

\newcommand{\wlr}[1]{\left( #1 \right)}

\newcommand{\wnode}[1]{\wfc{\ell}{#1}}

\newcommand{\wnorm}[1]{\| #1 \|}

\newcommand{\wref}[1]{(\ref{#1})}
\newcommand{\wrone}{\mathds R}

\newcommand{\wrm}[1]{\mathrm{#1}}

\newcommand{\wrounde}[1]{\wfc{\wrm{fl}}{#1}}

\newcommand{\wset}[1]{{\left\{ #1 \right\}}}

\newcommand{\wst}[1]{<\!\!{#1}\!\!>}
\newcommand{\wstx}[1]{<\!{#1}\!>}
\newcommand{\wstxx}[1]{<{#1}>}

\newcommand{\wvec}[1]{\mathbf{#1}}

\begin{document}

\title{The stability of barycentric interpolation at the Chebyshev points of the second kind}
\author{Walter F. Mascarenhas, \\ walter.mascarenhas@gmail.com\thanks{
Instituto de Matem\'{a}tica e Estat\'{i}stica, Universidade de S\~{a}o Paulo, 
Cidade Universit\'{a}ria, Rua do Mat\~{a}o 1010, S\~{a}o Paulo SP, Brazil. CEP 05508-090.
Tel.: +55-11-3091 5411, Fax: +55-11-3091 6134,walter.mascarenhas@gmail.com}           
}
\maketitle
\begin{abstract}
We present a new analysis of the stability of the first and second
barycentric formulae for interpolation at the Chebyshev points of the second kind.
Our theory shows that the second formula is more stable than previously thought and our
experiments confirm its stability in practice.
We also extend our current understanding regarding the accuracy problems of the first barycentric formula.
\end{abstract}

\section{Introduction}
 In 1972, Herbert Salzer first treated a special case of the remarkable {\it second barycentric formula} \cite{SALZER}:
\begin{equation}
\label{main}
\wfc{b_n}{t;f} :=
\frac{{\sum_{i = 0}^n} \frac{w_i f_i}{t - x_i}}{{\sum_{i =0}^n} \frac{w_i}{t - x_i}}
\end{equation}
for interpolating $f: [x_0,x_n] \mapsto \wrone{}$ in the nodes
$x_0 < x_1 < \dots < x_{n-1} < x_n$, with $f_i = \wfc{f}{x_i}$. This formula has its origins in Taylor \cite{TAYLOR},
which considered equally spaced $x_i$.
It leads to a function $b_n$ which is a polynomial in  $t$
when the weights $w_i$ are chosen as the numbers
\begin{equation}
\label{def_lambda}
\wfcc{\lambda_i}{\wvec{x}} := \prod_{j \neq i} \frac{1}{x_i - x_j}.
\end{equation}
We could also choose $w_i = \kappa_n \wfcc{\lambda_i}{\wvec{x}}$ for any constant $\kappa_n \neq 0$ independent
of $i$, because we do not change the right-hand side of \wref{main} when we multiply its numerator
and denominator by $\kappa_n$.
Salzer looked at the Chebyshev points of the second kind, which are given by
$x_i = - \wfc{\cos}{i \pi/n}$. Like Marcel Riesz in 1916 (see \cite{ATAP}, p.36), Salzer showed
that for these nodes we have
\begin{equation}
\label{riez}
\wfcc{\lambda_0}{\wvec{x}} = \frac{\wlr{-1}^n}{n}  2^{n-2}, \hspace{0.3cm}
\wfcc{\lambda_n}{\wvec{x}} = \frac{1}{n} 2^{n - 2} \hspace{0.3cm} \wrm{and} \hspace{0.3cm}
\wfcc{\lambda_i}{\wvec{x}} = \frac{\wlr{-1}^{n - i}}{n} 2^{n - 1} \ \ \wrm{for} \ \ 0 < i < n.
\end{equation}
As a result, when the $x_i$ are the Chebyshev points of the second kind we can write the second
formula concisely by dividing the $\lambda_i$ in \wref{def_lambda} by $\kappa_n = 2 \wfcc{\lambda_0}{\wvec{x}}$
and using the following {\it simplified weights} $\gamma_i$ as $w_i$ in \wref{main}:
\begin{equation}
\label{theoGamma}
\gamma_0 := 1 / 2, \hspace{0.5cm} \gamma_n := \wlr{-1}^n /2 \hspace{1.0cm} \wrm{and} \hspace{1.0cm}
\gamma_i := \wlr{-1}^{i}, \ \ \wrm{for} \ \ i = 1,\dots, n - 1.
\end{equation}

Taylor and Salzer studied the second barycentric formula
because of its simplicity and the accurate results it yielded.
In the four decades that have passed since the publication of Salzer's paper, the number of nodes $n$ considered
to be large has been raised from the hundreds to the million range.
Today, barycentric formulae are important tools in interpolation and also in spectral methods for the solution of partial differential equations.
Therefore, a deep understanding of
the numerical properties of formulae \wref{main}--\wref{theoGamma} is a
worthwhile scientific goal.

 In Salzer's time, errors of the order $n^2 \sup \wabs{f} \epsilon$ due to rounding, where $\epsilon \approx 10^{-16}$ is machine precision and we interpolate at $n + 1$ nodes, would have been considered satisfactory.
For one hundred nodes this corresponds to $10^{-12} \sup \wabs{f}$
and is reasonable.
Today however, for $n$ in the million
range considered in the fifth chapter of \cite{ATAP}, an error of $n^2 \sup \wabs{f} \epsilon$
corresponds to $10^{-4} \sup \wabs{f}$ and  is unacceptable.
 In this article, we show that proper implementations of the second barycentric formula
for $t \in [-1,1]$ with the simplified weights lead
to backward errors well below $n^2 \sup \wabs{f} \epsilon $.
We also present a bound on the forward error for functions with moderate derivatives.
We obtained these bounds by looking at the second formula from a new perspective.
Instead of writing it in the standard way, we use
an alternative formulation which improves
its stability when we use the simplified weights.

The main motivation for our analysis is the case discussed by Salzer, in
which the $x_i$ are the Chebyshev points of the second kind and
the rational function $\wfc{b_n}{t;f}$ is a polynomial.
However, in order to understand the numerical aspects of Salzer's case
we must first consider
the distinction between the abstract nodes $x_i$
and their rounded representation

\begin{equation}
\label{old_round}
\hat{x_i} = \wrm{a \ IEEE754 \ double \ precision \ number \ such \ that } \wabs{x_i - \hat{x}_i}
\wrm{\ is \ of \ order } \  10^{-16} \wabs{x_i}.
\end{equation}
Since the nodes are rounded in practice, our analysis of the numerical stability of the second barycentric
formula considers general families of sorted interpolation nodes
\begin{equation}
\label{sortedNodes}
-1 = \hat{x}_0 < \hat{x}_1 < \dots < \hat{x}_{n-1} < \hat{x}_n = 1.
\end{equation}
The rounded Chebyshev points of the second kind are not in
harmony with the simplified weights
and the corresponding second barycentric formula yields a rational function
which is not a polynomial. The consequences of this change in the approximation
of $f$ by $b_n$ are discussed in \cite{WERNER}
but they do not concern us in this paper since we focus solely on the numerical aspects of the resulting formulae.

The literature does not pay due attention to the case of simplified weights in \wref{theoGamma} and $t \in [-1,1]$ that we consider.
For instance, \cite{WEBBX} is concerned with extrapolation and not interpolation whereas \cite{HIGHAM} considers  generic nodes.
Unfortunately, as we show in section \ref{secFirst} below,
 the situation is less favourable for the first barycentric formula.
In  section \ref{secFirst}, we report significant problems with this formula.
We show that implementing the first formula with the Chebyshev points of the second
kind rounded as usual may lead to  errors of order
$n^2 \sup \wabs{f} \epsilon$ for $\wfc{f}{t} = \wfc{\sin}{t}$ and $t$ near the nodes.
For the sake of brevity and clarity, we focus on the big picture regarding numerical stability for interpolation (providing the technicalities in the appendices)
and we eschew any specifics regarding the efficient implementation of our version of the second formula \wref{main}. However, it is possible to implement the first
barycentric formula efficiently and accurately, but this implementation requires that we handle the
nodes in an unusual way. In the following sections,
we start with general remarks and present experiments to illustrate their point. Next, we informally explain the mathematics underlying the experiments, and we state
a lemma or a theorem if applicable. Section \ref{sec_the_end} summarizes our findings and discusses future work. 

\section{The stability of the second formula}
\label{secSecond}
We reformulated the second barycentric formula for the simplified weights in  \wref{theoGamma}, involving a series of mathematical expressions detailed in section \ref{subOverview}.
 Here, we explore the consequences of our formulation, mainly its property of backward stability with respect to perturbations in the function values and the resulting small forward errors for functions with moderate derivatives.

The abstract in \cite{HIGHAM}
mentions that the error analysis
for the second formula is less favorable than the one for the first formula, but we
should not infer from this abstract that the second formula is less stable than the first one.
Higham provides only upper bounds on the error.
Sometimes error bounds are realistic, at other times they are not.

Our theory contains two theorems. In order to motivate them, we start with
an experiment that illustrates the numerical stability of the second formula.
It considers the highly oscillating function
$\wfc{f}{t} = \wfc{\sin}{10^5 t}$, $f_i = \wfc{f}{x_i}$, the Chebyshev points of the second kind
and the simplified weights. In this case
the function $\wfc{b_n}{t;f}$ defined by \wref{main}
is a $n$-th degree polynomial in $t$ which interpolates $f$
at the nodes $x_i$. Therefore, we can estimate the error in approximating
$\wfc{f}{t}$ by $\wfc{b_n}{t;f}$ by (see pg. 187 of \cite{HENRICI}):
\begin{equation}
\label{bound5}
\wabs{\wfc{b_n}{t;f} - \wfc{f}{t}} \leq \frac{\max_{-1 \leq \xi \leq 1}\wabs{\wfc{f^{\wlr{n+1}}}{\xi}}}{\wlr{n+1}!}
\wabs{\wnode{t}},
\end{equation}
for $\wnode{t} = \prod_{i = 0}^n \wlr{t - x_i}$.
Using the chain rule, Stirling's formula and the bound $\wabs{\wnode{t}} \leq 2^{1-n}$ presented
in the first page of \cite{SALZER}, we obtain this back-of-the-envelope simplification of \wref{bound5}:
\begin{equation}
\label{bound6}
\wabs{\wfc{b_n}{t;f} - \wfc{f}{t}} \preceq \wlr{ \frac{10^5 e}{2 n} }^n.
\end{equation}
This heuristic bound does not say much for $n$ smaller than $10^5 e / 2$. However,
if $n$ is a few percent larger than $10^5 e / 2$ then it tells us
that the error is well below the precision of our computers.
For $n$ about one million, the bound \wref{bound6} shows that
$b_n$ and $f$ are indistinguishable.

The bound \wref{bound5} shows how accurately
$b_n$ approximates $f$.
Here,
we are concerned with the {\it stable} evaluation of $b_n$,
so that this accuracy is not ruined by rounding errors. Table \ref{tableErrorSin}
illustrates the numerical stability of $b_n$\footnote{The experimental details are described in appendix \ref{apExperiments}.}.
It compares $b_n$'s results with those obtained using
the following naive implementation of $f$ in  \verb C++ :
\begin{verbatim}
double sin_10To5_X(double t) { return sin( 100000 * t ); }
\end{verbatim}
The type \verb double  in \verb C/C++  implements
the IEEE754 double precision, which yields a machine precision $\epsilon \approx {10}^{-16}$,
 which is also used by Matlab and corresponds to \verb real*8  in \verb Fortran .

{\small
\begin{table}[ht]
\centering
\caption{Error statistics for the evaluation of $\wfc{\sin}{10^5 t}$ in $10^6$ equally spaced points from $t = -1$ to
$t = 1$.
}\begin{tabular}{r|cc|cc}
\hline\hline\\[-0.25cm]
Number           & \multicolumn{2}{|c|}{C++ function}               & \multicolumn{2}{|c}{Second barycentric formula $b_n$}      \\
of nodes         & max                   & mean$\pm$std. dev.           & max              & mean$\pm$std. dev.     \\
\hline\\[-0.27cm]
$10^5 + 1$           & $7.3\times10^{-12}$ & $1.3\pm1.5\times 10^{-12}$ & $0.67$                & $0.02\pm0.07$              \\
$5 \times 10^5 + 1$  & $7.3\times10^{-12}$ & $1.3\pm1.5\times 10^{-12}$ & $6.0 \times 10^{-12}$ & $4.6\pm5.6\times 10^{-13}$ \\
$10^6 + 1$           & $7.3\times10^{-12}$ & $1.3\pm1.5\times 10^{-12}$ & $3.7 \times 10^{-12}$ & $3.3\pm3.9\times 10^{-13}$ \\
$25 \times 10^6 + 1$ & $7.3\times10^{-12}$ & $1.3\pm1.5\times 10^{-12}$ & $7.0 \times 10^{-13}$ & $6.5\pm7.7\times 10^{-13}$ \\
\end{tabular}
\label{tableErrorSin}
\end{table}
}

The rows related to the \verb C++  function in Table \ref{tableErrorSin} are equal.
They correspond to the same instances of $t$ and the same function. The third and fourth
rows in the columns corresponding to $b_n$ are quite different, as one would
expect from \wref{bound6}. The last rows in Table \ref{tableErrorSin}
show that not only is the second barycentric formula
as accurate as the \verb C++  function, it is even more accurate.

\subsection{Formal results regarding the second formula}
\label{formalSecond}
This section analyzes the numerical stability of the
second barycentric formula with the
simplified weights in \wref{theoGamma} for $t \in [-1,1]$. We make the usual assumptions
about floating point arithmetic described in \cite{HIGHAM_BOOK}.
We present two theorems: the first one expresses the backward error in terms of relative perturbations in the function values, in the spirit of \cite{HIGHAM} and the second theorem is useful for functions with moderate derivatives. The second theorem is more powerful than the first because it provides a small forward error bound, without a quadratic dependency on $n$.

We assume that we can sum numbers $a_0,\dots, a_n$
to obtain a numerical sum $\hat{\sum}_{i = 0}^n a_i$ with
\begin{equation}
\label{sumA}
\hat{\sum}_{i = 0}^n a_i = \sum_{i = 0}^n \wlr{1 + \delta_i \sigma_n \epsilon } a_i,
\end{equation}
where $\epsilon$ is machine precision, $\wabs{\delta_i} \leq 1$ and $\sigma_n$
may depend on $n$ but is unrelated to the $a_i$. Several strategies yield
reasonably small $\sigma_n$ in \wref{sumA}, including
the naive one, which leads to $\sigma_n = n$ for the $n$ and $\epsilon$ relevant in practice. Readers who need extra accuracy
could consider Kahan's summation \cite{KAHAN},
which yields $\sigma_n = 2 + \wfc{O}{n \epsilon}$,
or the methods by Rump or Priest
\cite{RUMP,PRIESTA,PRIESTB}.
We simply assume that the readers can sum according to \wref{sumA}
and let them use their own $\sigma_n$ in the following theorems (we also assume that $1 \leq \sigma_{n} \leq \sigma_{n+1})$:
\begin{thm}
\label{thmBackward}
Consider $n \geq 2$, nodes $\hat{x}_i$ as in \wref{sortedNodes} and define $m := \wfloor{(n+1)/2}$.
Assume we can sum as in  \wref{sumA}, with $\sigma_n$
and $\sigma_m$ such that $\wlr{8 + 1.01 \sigma_n + 1.01 \sigma_m} \epsilon < 0.01$.
If $2 + \hat{x}_1$ and $2 - \hat{x}_{n-1}$ are floating point numbers then we
can evaluate the second barycentric formula in \wref{main} with
the simplified weights in \wref{theoGamma} so that the computed value
$\wfc{\hat{b}_n}{t;f}$ is equal to $\wfc{b_n}{t;\tilde{f}}$ for a vector
$\tilde{f}$  with
\begin{equation}
\label{backBound}
\wabs{\tilde{f}_i - f_i} \leq \wlr{8.1 + 1.1 \sigma_n + 1.1 \sigma_m} \wabs{f_i} \epsilon.
\wend{}
\end{equation}
\end{thm}

\begin{thm}
\label{thmMain}
Consider $n \geq 2$, nodes $\hat{x}_i$ as in \wref{sortedNodes} and define $m := \wfloor{(n+1)/2}$.
Assume that we can sum as in \wref{sumA}, with $\sigma_m$ such that $\wlr{20 + 2.02 \sigma_m} \epsilon < 0.01$.
If $2 + \hat{x}_1$, $2 - \hat{x}_{n-1}$ and $\hat{x}_i + \hat{x}_{i - 1}$, for $i = 1,\dots,n- 1$,
are floating point numbers then we can evaluate the second barycentric formula \wref{main} with
the simplified weights in \wref{theoGamma}
so that the computed value $\wfc{\hat{b}_n}{t;f}$ is equal to $\wfc{b_n}{t;\hat{f}}$ for a vector
$\hat{f}$ with \footnote{$\wnorm{v}_\infty = \max_{i} \wabs{v_i}$ is the sup norm of the vector $v$ and,
by convention, $f_{-1} = f_{n+1} = 0$.}
\begin{equation}
\label{backBoundMain}
\wabs{\hat{f}_i - f_i} \leq \wlr{16.6 + 2.1 \sigma_m} \max \wset{ \wabs{f_{i-1}} + \wabs{f_i}, \wabs{f_i} + \wabs{f_{i + 1}} } \epsilon.
\end{equation}
Moreover,
\begin{equation}
\label{forwardBound}
\wabs{\wfc{\hat{b}_n}{t;f} - \wfc{b_n}{t,f}} \leq
 \wlr{45.6 \wnorm{f}_\infty +  37.4 \wnorm{\Delta  \! f }_\infty + 6.2 \wnorm{f}_\infty \sigma_m +
          4.1 \wnorm{\Delta \! f }_\infty \sigma_m}
 \epsilon,
\end{equation}
where $\Delta \! f$ is the vector with entries
$\Delta \! f_i = \wlr{f_{i} - f_{i-1}}/\wlr{x_{i} - x_{i-1}}$. \wend{}
\end{thm}

To appreciate these theorems, the reader should compare the bounds they provide with the ones
presented in \cite{HIGHAM} for general weights.
Our bounds are smaller,
even if we take into account that they may hide a factor of $n$ in
$\sigma_n$ and $\sigma_m$. However, this is not the most important difference
between our results and \cite{HIGHAM}. What matters most is
our proof of backward stability with respect to the function values. Higham's
work does not lead to this conclusion,
but does not exclude it either.
Our work is more informative regarding this particular case, whereas Higham's work is more appropriate within the general setting.

The forward bound \wref{forwardBound} is
small if the $f_i$ come from a function $f$ with a moderate derivative, in view of
Lemma \ref{lemSmallDf} below. Combining this lemma with Theorem
\ref{thmMain}, we conclude that if we use Kahan's summation and compute
the $f_i$  with a small error, then we obtain
a result with an error of the order of a small constant times $\epsilon$ for all
$n$ relevant in practice.

\begin{lem}
\label{lemSmallDf} If the function $f: [-1,1] \to \wrone{}$ is
differentiable and $f_i = \wfc{f}{x_i} + \delta_i$ then
\[
\wnorm{\Delta \! f}_\infty \leq \wnorm{f'}_\infty + \frac{2 \wnorm{\delta}_\infty}{\min_{j = 1, n} \wlr{x_{j} - x_{j-1}}}.
\]
\wend{}
\end{lem}

The hypothesis, that $2 + \hat{x}_1$, $\hat{x}_{i} + \hat{x}_{i-1}$ and
$2 - \hat{x}_{n-1}$ are floating point numbers in Theorem \ref{thmBackward}
and Theorem \ref{thmMain}, is inconvenient. Unfortunately, we did not
find a way to replace it with anything simpler without introducing
terms that depend on the node spacing in our bounds.
However,
it is easy to obtain rounded nodes satisfying the hypothesis of our theorems.
In appendix \ref{apRoundingNodes}, we present \verb C++  code
that produces rounded nodes $\hat{x}_i$ that satisfy the
hypothesis of our theorems
 for which we can prove, under reasonable assumptions, that
$\wabs{x_i - \hat{x}_i} \leq 2.53 \wabs{x_i} \epsilon$
for the nodes $x_i = - \wfc{\cos}{i \pi / n}$ with $n \leq 10^9$.
In fact, our experiments indicate that our rounded nodes satisfy $\wabs{x_i - \hat{x}_i} \leq 2 \wabs{x_i} \epsilon$
in this case.
The same experiments show that nodes obtained from the
\verb C++   evaluation of such $x_i$ using the more accurate formula $x_i = \wfc{\sin}{(2 i - n) \pi / 2 n}$
may have errors of size $1.6 \wabs{x_i}  \epsilon$. Therefore, our rounded nodes
are almost as accurate as the ones used today.

\section{The first barycentric formula is problematic}
\label{secFirst}
The first barycentric formula can be written as:
\begin{equation}
\label{first}
\wfc{a}{t} :=  \frac{\wlr{-1}^n}{n}  \ \times \ \wlr{2^{n-1}
\prod_{i = 0}^n \wlr{t - x_i} } \ \times \
\sum_{i = 0}^n \frac{w_i f_i}{t - x_i}.
\end{equation}
With appropriate weights $w_i$, this formula is equivalent to
the second one in exact arithmetic. However, they are quite different from the numerical point of view.
We consider the cases in which the $w_i$ are the simplified weights in \wref{theoGamma} and in which the $w_i$ are obtained by evaluating the following quantity  numerically:
\begin{equation}
\label{firstNu}
\nu_i := \wlr{-1}^n \frac{n}{2^{n-1} \prod_{j \neq i} \wlr{\hat{x}_i - \hat{x}_j}}.
\end{equation}
These $\nu_i$ are the weights that turn the right-hand side of
\wref{first} into the usual first formula when we replace the exact arithmetic nodes $x_i$
by the rounded nodes $\hat{x}_i$.
 We discuss the following issues:
\begin{enumerate}
\item Overflow and underflow.
\item The instability of the first formula evaluated with the simplified weights in \wref{theoGamma}.
\item The inconvenience of the first formula with the weights $w_i = \nu_i$ in \wref{firstNu}.
\end{enumerate}

 In subsection \ref{subFirst}, we show that if we use Chebyshev points of the second kind rounded as in \wref{old_round}, then this formula is inaccurate for $t$ near the nodes.
In this case the errors can be of order $n^2 \wnorm{f}_\infty \epsilon$,
which is unacceptable for large $n$. Due to these large errors, we did consider using the
 numerically computed weights by evaluating \wref{firstNu}.

Indeed, our experiments show that the weights in \wref{firstNu} lead
to significantly more accurate results.
However, their use leads to performance issues which we discuss in section \ref{subTheLastOne}.  The overall conclusion is that,
for stability and performance reasons, we should prefer the second barycentric
formula to the first one for interpolation in the
Chebyshev points of the second kind rounded as usual.

\subsection{The triviality of overflow and underflow in the first barycentric formula}
\label{subOverflow}
Concerns and strategies regarding the overflow and underflow in the intermediate steps of the computation of the product $\prod_{i = 0}^n \wlr{t - x_i}$
in the evaluation of the first barycentric formula are outlined on page 509 of \cite{BERRUT} and in the fifth chapter of \cite{ATAP}.
 This section and appendix \ref{apOverflow} present our strategy that we consider to be simpler and more accurate.

As in \wref{first}, instead of using only the product $\prod_{i = 0}^n \wlr{t - x_i}$, we consider
\begin{equation}
\label{firstProduct}
2^{n-1} \prod_{i = 0}^n \wlr{t - x_i},
\end{equation}
changing the factor $2^{n-1}$ in this formula if our interpolation interval is not $[-1,1]$.
From this perspective, the cure to overflow and underflow is indeed so simple that we do not consider
them to pose real problems with the first formula.

Our solution efficiently scales the intermediate products by appropriate powers of two, keeping track of this scaling and without introducing any rounding errors.
Our implementation is written in \verb C++ , but it can be adapted to similar languages or Matlab.
The solution is based on two Ansi \verb C  functions called
\verb frexp  and \verb scalbln  declared as follows:
\begin{verbatim}
 double frexp(t, int* exp)
 double scalbln(double v, long int exp).
\end{verbatim}
\noindent
\verb scalbln  scales its argument \verb t  by \verb 2^exp .
If $t \neq 0$  then \verb frexp  returns \verb y $ \in [0.5,1) \cup (-1,-0.5]$
and sets the exponent \verb *exp  so that \verb t \verb = \verb y \verb 2^*exp .
If \verb t  is zero then it returns zero and sets \verb *exp  to zero.
Appendix \ref{apOverflow}  contains the \verb C++  implementation
of our solution.
We expect it to be as efficient as the naive implementation of \wref{firstProduct} due to our parsimonious scaling.
We now present an experiment comparing the speed and
accuracy of our solution, which we call \verb Scaling ,
with the following alternatives
\footnote{The experimental settings are described in appendix \ref{apExperiments}}:
\begin{enumerate}

\item \verb Naive  implements the product as a simple loop, ignoring overflow and underflow.
\item \verb LogSum  computes \wref{firstProduct} by taking the sums of the logs of the factors and then
exponentiating the resulting sum, as proposed in the fifth chapter of \cite{ATAP}.
\item \verb IppLogSum  optimizes \verb LogSum  by using
Intel's Integrated Performance Primitives.
\item \verb Grouped  \verb logs  is similar to \verb Scaling , with calls to \verb frexp  and \verb scalbln  replaced
by \verb log  and \verb exp . As \verb Scaling , \verb Grouped  \verb logs  amortizes the cost of the logs by
computing them only once per twenty products.
\end{enumerate}

We present the experimental results in three tables.
 The first concerns performance whereas the other two are about accuracy.
The performance results seem
to depend on our experimental settings and should only be taken as an indication of efficiency.
The accuracy results do not depend on our settings and our
experiments present evidence that \verb Scaling  is more accurate than sums of logs.
 In fact, for large $n$, the experiment indicates that summing logs leads to inaccurate results and illustrates the poor performance of the \verb Naive  approach.

{\small
\begin{table}[ht]
\caption{Relative time to compute $10^5$ products}
\centering
\begin{tabular}{c|ccccc}
\hline\hline\\[-0.23cm]
Number  &             &                &              &    & \\
of nodes     & \verb Naive & \verb Scaling  & \verb LogSum & \verb IppLogSum & \verb Grouped  \verb logs  \\
\hline\\[-0.25cm]
$10^3$    & 1            & 0.34           &  6.0         &  1.1 &  0.84 \\
$10^4$    & 1            & 0.84           &  19.8        &  2.9 &  1.84 \\
$10^5$    & 1            & 0.41           &  8.5         &  1.6 &  0.83 \\
$10^6$    & 1            & 0.41           &  7.3         &  2.3 &  0.76 \\
\end{tabular}
\label{tableTimings}
\end{table}
}

Table \ref{tableTimings} expresses times as multiples of the time taken by \verb Naive  and shows that \verb Scaling  can be much faster than \verb Naive .
We did not expect this and do not consider it to be evidence that \verb Scaling 's
performance is superior to \verb Naive 's in any other combination of programmer,
software and hardware.
Performance in today's computers depends
on several factors and the use of the appropriate libraries to take advantage of
them. This  can be seen by comparing columns \verb LogSum  and
\verb IppLogSum  in Table \ref{tableTimings}.
We tried to be fair with all methods, using
an optimized library 
 to implement \verb IppLogSum  and implementing the
\verb Grouped  \verb logs  strategy in a way that differs from our own
strategy only by the use of \verb log/exp  instead of \verb frexp/scalbln .
However, we must recall that our performance results depend on our particular settings.

Tables \ref{tableMaxError} and \ref{tableMeanError} below present statistics
on the errors in the evaluation of products with factors given by the term inside
parentheses in equation \wref{first} for $10^5$ $t$ chosen randomly
in $[-1,1]$.
They show that \verb Scaling  handles overflow and underflow properly
and is accurate, with the rounding errors increasing linearly with
the number of nodes. The alternatives do not perform so well: the naive
method is unacceptable and the strategies using sums of logs are significantly
less accurate.

{\small
\begin{table}[ht]
\caption{Maximum relative error in $10^5$ products}
\centering
\begin{tabular}{c|cccccc}
\hline\hline\\[-0.23cm]
Number   &                 &                    &               &   & \\
of nodes &   \verb Naive   & \verb Scaling      & \verb LogSum  & \verb IppLogSum & \verb Grouped  \verb logs     \\
\hline\\[-0.25cm]
$10^3$   &  $490$         &  $2.5 \times 10^{-14}$   & $3.8 \times 10^{-12}$  & $1.0 \times 10^{-12}$  & $7.8 \times 10^{-15}$ \\
$10^4$   & $\infty$       &  $2.1 \times 10^{-13}$   & $1.7 \times 10^{-10}$  & $2.4 \times 10^{-11}$  & $1.7 \times 10^{-11}$ \\
$10^5$   & $\infty$       &  $2.0 \times 10^{-12}$   & $4.3 \times 10^{-09}$  & $6.7 \times 10^{-10}$  & $5.0 \times 10^{-10}$ \\
$10^6$   & $\infty$       &  $2.0 \times 10^{-11}$   & $1.3 \times 10^{-07}$  & $1.7 \times 10^{-08}$  & $1.6 \times 10^{-08}$ \\
\end{tabular}
\label{tableMaxError}
\end{table}
}

{
\begin{table}[ht]
\caption{Average relative error $\pm$ standard deviation in $10^5$ products}
\centering
\tabcolsep=0.11cm
\begin{tabular}{c|ccccc}
\hline\hline\\[-0.23cm]
Number   &                 &                            &                            &  & \\
 of nodes &  \verb Naive   & \verb Scaling              & \verb LogSum               & \verb IppLogSum & \verb Grouped  \verb logs \\
\hline\\[-0.25cm]
$10^3$    & $ 0.3\pm1.4 $  & $2.8\pm4.2\times 10^{-15}$ & $5.4\pm4.5\times 10^{-13}$ &  $5.1\pm1.1\times 10^{-13}$ &  $1.3\pm1.0\times 10^{-13}$   \\
$10^4$    & $ \infty $     & $2.1\pm4.3\times 10^{-14}$ & $1.6\pm1.5\times 10^{-11}$ &  $5.3\pm2.6\times 10^{-12}$ &  $2.6\pm2.0\times 10^{-12}$   \\
$10^5$    & $ \infty $     & $1.8\pm4.3\times 10^{-13}$ & $6.1\pm5.5\times 10^{-10}$ &  $8.9\pm7.3\times 10^{-11}$ &  $8.2\pm6.5\times 10^{-11}$  \\
$10^6$    & $ \infty $     & $1.7\pm4.4\times 10^{-12}$ & $1.8\pm1.5\times 10^{-08}$ &  $2.2\pm1.9\times 10^{-09}$ &  $2.5\pm2.0\times 10^{-09}$   \\
\end{tabular}
\label{tableMeanError}
\end{table}
}

\subsection{Problems with the first formula}
\label{subFirst}
We now explain why we believe the simplified weights in \wref{theoGamma} should not be used
in combination with Chebyshev points of the second kind rounded as in \wref{old_round} when evaluating
the first barycentric formula \wref{first} with large $n$.
We show that, in this case, the first formula yields results of widely varying accuracy.
The results will all look deceptively nice in plots and we will not realize that some of them have
errors worse by orders of magnitude than could be incurred using the second formula.
In the end, this may be even less desirable than having numbers that are evidently wrong.

Our experiment is designed to test the stability of the first formula for $t$ near the nodes when we use the simplified weights and Chebyshev points of the second kind rounded as in \wref{old_round}.
The experiment compares three approaches. The first one implements the first barycentric
formula \wref{first} with the simplified weights in \wref{theoGamma} and nodes $x_i$ obtained
by rounding $\wfc{\sin}{\wlr{2 i - n} \pi / {2n}}$, which is a more accurate form of
writing $- \wfc{\cos}{i \pi/n}$. The second implements the same formula \wref{first}
with the same nodes, but instead of the simplified weights we use the ones
obtained by evaluating $\nu_i$ in \wref{firstNu} numerically.
The third approach is based on the results in section \ref{secSecond}.
We call it Stable.
The points at which we evaluate the formulae are the key aspect  of this experiment
and  are in fact extremely close to the nodes.
We take the hundred nodes  $\hat{x}_{n-100}$, $\dots$, $\hat{x}_{n-1}$ and for each
of them we consider the nearest $10^4$ floating point numbers, to the left and to the right.
The results are presented in  Table \ref{tableNever}, providing strong evidence for the instability of the first formula with simplified weights and suggesting a quadratic dependency of the error on $n$.

Using the better weights in \wref{firstNu}, we find our results are much improved but still inferior to those
obtained with the second formula.

{\small
\begin{table}[ht]
\caption{Maximum error for $2 \times 10^6$ points near the nodes for $\wfc{f}{t} = \wfc{\sin}{t}$.}
\centering
\begin{tabular}{c|ccc}
\hline\hline\\[-0.25cm]
Number of       & \multicolumn{3}{c}{Formula} \\
\cline{2-4}\\[-0.25cm]
 nodes          & First, simplified weights & First, weights as in \wref{firstNu} &  Stable \\
\hline\\[-0.25cm]
$ 10^3 + 1$     & $7.3\times 10^{-12}$       & $9.7\times 10^{-15}$         & $4.0\times 10^{-16}$ \\
$ 10^4 + 1$     & $5.2\times 10^{-10} $      & $3.0\times 10^{-14}$         & $4.3\times 10^{-16}$ \\
$ 10^5 + 1$     & $8.7\times 10^{-08} $      & $1.2\times 10^{-13}$         & $4.2\times 10^{-16}$ \\
$ 10^6 + 1$     & $6.0\times 10^{-06} $      & $3.4\times 10^{-13}$         & $4.8\times 10^{-16}$
\end{tabular}
\label{tableNever}
\end{table}
}

{\small
\begin{table}[ht]
\caption{Error $\pm$ standard deviation for $2 \times 10^6$ points near the nodes for $\wfc{f}{t} = \wfc{\sin}{t}$.}
\centering
\begin{tabular}{c|ccc}
\hline\hline\\[-0.25cm]
Number of       & \multicolumn{3}{c}{Formula} \\
\cline{2-4}\\[-0.25cm]
 nodes          & First, simplified weights & First, weights as in \wref{firstNu} &  Stable \\
\hline\\[-0.25cm]
$ 10^3 + 1$     & $2.3\pm6.2\times 10^{-13}$      & $1.6\pm1.2\times 10^{-15}$         & $7.2\pm5.5\times 10^{-17}$ \\
$ 10^4 + 1$     & $1.7\pm4.8\times 10^{-11} $     & $5.2\pm3.9\times 10^{-15}$         & $7.3\pm5.5\times 10^{-17}$ \\
$ 10^5 + 1$     & $2.5\pm7.9\times 10^{-09} $      & $2.3\pm1.7\times 10^{-14}$         & $7.4\pm5.6\times 10^{-17}$ \\
$ 10^6 + 1$     & $1.7\pm4.8\times 10^{-07} $      & $6.2\pm4.6\times 10^{-14}$         & $7.2\pm5.5\times 10^{-17}$
\end{tabular}
\label{tableNeverStd}
\end{table}
}

 In order to provide a theoretical explanation of the experimental data, we consider the effects of rounding in the first formula.
The main reference on this subject is Higham's article \cite{HIGHAM}
\footnote{We also recommend his book \cite{HIGHAM_BOOK} for a broader
view of the issues we discuss here.}.
Higham's definition  of backward stability is based on relative perturbations in the function values.
 This definition is reasonable but we outline some of its limitations.
For example, it says nothing about barycentric interpolation in regard
to the definition
of  ``{\it backward stability}'' according to which an algorithm
to evaluate a function $f$ is ``{\it backward stable}'' if the
computed value of $\wfc{f}{x}$ is the exact
value of $\wfc{f}{y}$ for some $y$ near $x$. The difference
between backward stability in  Higham's text and backward stability
in the latter sense is evident when $x$ is a global
maximizer of the function $f$ we want to approximate. In this case,
there may not exist a $y$ that satisfies the requirements of backward
stability in the latter sense, but we could fulfill the requirements of backward stability
in  Higham's sense by perturbing $f$.

We now extend Higham's work by considering the effects of rounding the nodes in
the first barycentric formula. We claim that these errors indirectly affect the weights
and that by doing so they cause instability. In fact, equation (3.2) in \cite{HIGHAM}
shows that Higham's version of the first formula is equivalent to our equation
\wref{first} with weights $\nu_i$ given by equation \wref{firstNu}, since he assumes that the nodes
defining the weights are the ones used in the computation and these nodes
are the $\hat{x}_i$. The simplified weights $\gamma_i$ in \wref{theoGamma}
are slightly different from the $\nu_i$, in that the
$\gamma_i$ correspond to exact nodes. What are the
consequences of this mismatch between $\gamma_i$ and $\nu_i$?
A naive answer to this question would be ``It does not matter; the effect of rounding in the nodes
is of the order machine epsilon and therefore negligible.'' However, these rounding errors
{\it do} matter
for large $n$, because $\gamma_i - \nu_i$ may be of order $n^2 \epsilon$ for the
Chebyshev points of the second kind rounded as in \wref{old_round}.

 The large errors in the first column of Table \ref{tableNever} occur because we actually evaluate
\begin{equation}
\label{aTilde}
\wfc{\tilde{a}}{t, \hat{\wvec{x}}} := \wlr{-1}^n \frac{2^{n-1}}{n}
\wlr{ \prod_{i = 0}^n \wlr{t - \hat{x}_i} }
\sum_{i = 0}^n \frac{\gamma_i f_i}{t - \hat{x}_i},
\end{equation}
 instead of the first barycentric formula $\wfc{a}{t}$ in \wref{first}.
(We use a bold $\wvec{z}$ to indicate the vector $(z_0,z_1,\dots,z_n)^t$ and ignore rounding in $f_i$.)
As we now explain, it is likely that $\wfc{\tilde{a}}{t, \hat{\wvec{x}}}$  differs from
$\wfc{a}{t}$ by $\wfc{\Theta}{n^2 \epsilon}$ even in exact arithmetic
\footnote{As in \cite{CORMEN}, we use $\Theta$ to denote ``of the same order as'' and use $O$ to denote
``of order up to''.}.
We measure rounding errors and the distance of $t$ to $x_k$ in terms of
\begin{equation}
\label{alphaDelta}
\wfc{\alpha_{n,k}}{\hat{\vec{x}}}  :=  \sum_{i \neq k} \wabs{\frac{\hat{x}_i - x_i}{x_i - x_k}}
\hspace{0.7cm} \wrm{and} \hspace{0.7cm}
\wfc{\delta_{n,k}}{t}  :=  \sum_{i \neq k} \wabs{\frac{t - x_k}{x_i - x_k}}.
\end{equation}
Note
that $\wfc{\alpha_{n,k}}{\hat{\wvec{x}}}$ and $\wfc{\delta_{n,k}}{\hat{x}_k}$
converge to zero as the rounding errors $\wabs{\hat{x}_i - x_i}$ become very small and $\wfc{\delta_{n,k}}{t}$
converges to zero when $t$ approaches $x_k$.
 In the following arguments, the reader
should think of $\wfc{\alpha_{n,k}}{\hat{\wvec{x}}}$, $\wfc{\delta_{n,k}}{\hat{x}_k}$ and
$\wfc{\delta_{n,k}}{t}$ as very small numbers.

The next lemma relates the $\wfc{a}{t}$ we would like to compute
to the value $\wfc{\tilde{a}}{t,\hat{\wvec{x}}}$
which
we evaluate numerically when we use the first barycentric
formula \wref{first} with the Chebyshev points of the second kind for $t$ very close to $x_k$:

\begin{lem}
\label{lemChi} Consider $1 \leq k < n$. If $f_0,\dots,f_n$
 and $t \in [-1,1]$ are such that
\begin{equation}
\label{chiH}
\wabs{f_k} \geq \wnorm{f}_\infty/2,
\hspace{0.4cm}
\wfc{\alpha_{n,k}}{\hat{\wvec{x}}} \leq 1/{24},
 \hspace{0.4cm}
\wfc{\delta_{n,k}}{\hat{x}_k} \leq 1/{24}
 \hspace{0.4cm} \wrm{and} \hspace{0.4cm}
\wfc{\delta_{n,k}}{t} \leq 1/{24},
\end{equation}
for $x_i = - \wfc{\cos}{i \pi / n}$, then there exist $\wfc{\beta_{n,k,i}}{t,\hat{\wvec{x}}}$ and
$\wfc{\kappa_{n,k,i}}{t,\hat{\wvec{x}}}$ such that
\begin{equation}
\label{eqTA}
\wfc{\tilde{a}}{t,\hat{\wvec{x}}} = \wfc{a}{t} e^{-\wfc{r_{n,k}}{t,\hat{\wvec{x}}}},
\end{equation}
 for $\wfc{r_{n,k}}{t,\hat{\wvec{x}}}$  given by
\begin{equation}
\label{rnk}
 \sum_{i\neq k} \frac{\wfc{\beta_{n,k,i}}{t,\hat{\wvec{x}}} \wlr{\hat{x}_i - x_i} +
                 \wfc{\kappa_{n,k,i}}{t,\hat{\wvec{x}}} \wlr{\hat{x}_k - x_k}}{x_k - x_i}
 - \sum_{i\neq k} \frac{\wfc{\kappa_{n,k,i}}{t,\hat{\wvec{x}}}
                      \wlr{t - x_k} \wlr{\hat{x}_i - x_i}}{\wlr{x_k - x_i}\wlr{t - x_i}}
\end{equation}
and the functions $\beta_{n,k,i}$ and $\kappa_{n,k,i}$ are almost constant:
 \begin{equation}
 \label{dBeta}
 \wabs{ \wfc{\beta_{n,k,i}}{t,\hat{\wvec{x}}} - 1} \leq 1.05 \wfc{\delta_{n,k}}{t} +
 0.6 \wfc{\alpha_{n,k}}{\hat{\wvec{x}}}
  \end{equation}
 \begin{equation}
 \label{dKappa}
 \wabs{\wfc{\kappa_{n,k,i}}{t,\hat{\wvec{x}}} - \frac{\gamma_i f_i}{\gamma_k f_k}}
  \leq
   2.4 \wfc{\alpha_{n,k}}{\hat{\wvec{x}}}+ 7.2 \wfc{\delta_{n,k}}{t} + 4 \wfc{\delta_{n,k}}{\hat{x}_k} +
   4 \wfc{\alpha_{n,k}}{\hat{\wvec{x}}} \wfc{\delta_{n,k}}{t}
  .
\wend{}
 \end{equation}
\end{lem}

When $f_k$ is not small, if the rounding errors are very small and if $t$ is very close to $x_k$
then $\wfc{\alpha_{n,k}}{\hat{\wvec{x}}}$, $\wfc{\delta_{n,k}}{t}$ and
$\wfc{\delta_{n,k}}{\hat{x}_k}$ are very small and Lemma \ref{lemChi},
through equations \wref{dBeta} and $\wref{dKappa}$, allows us to replace $\wfc{\beta_{n,k,i}}{t,\hat{\wvec{x}}}$ by $1$ and
$\wfc{\kappa_{n,k,i}}{t,\hat{\wvec{x}}}$ by $\gamma_i f_i / {\gamma_k f_k}$ in
\wref{rnk}.
We can then neglect the second-order terms $\wlr{t - x_k} \wlr{\hat{x}_i - x_i}$,
and write
\begin{equation}
\label{snk}
\wfc{r_{n,k}}{t,\hat{\wvec{x}}} \approx
\wfc{r_{n,k}}{x_k,\hat{\wvec{x}}} =
\wfc{s_{n,k}}{\hat{\wvec{x}}} =
\frac{1}{\gamma_k f_k}
\sum_{i\neq k} \frac{\gamma_k f_k  \wlr{\hat{x}_i - x_i} + \gamma_i f_i \wlr{\hat{x}_k - x_k}}{x_k - x_i}
\end{equation}
and estimate the relative error caused by rounding in the evaluation of $\wfc{a}{t}$,
because in these circumstances \wref{eqTA} and \wref{snk} show that, to leading order, this error is equal to
$\wfc{s_{n,k}}{\hat{\wvec{x}}}$.

If we think of the rounding errors $\hat{x}_i - x_i$ as independent random variables with mean zero and
standard deviation
$\sigma_{n,i} \approx \wabs{x_i} \epsilon$, then  $\wfc{s_{n,k}}{\hat{\wvec{x}}}$ is a random variable with mean
zero and standard deviation
\begin{equation}
\label{var}
\wfc{\sigma}{\wfc{s_{n,k}}{\hat{\wvec{x}}}} \approx \epsilon \sqrt{ \sum_{i \neq k} \wlr{\frac{x_i}{x_k - x_i}}^2  +
\wlr{\sum_{i \neq k} \frac{\gamma_i f_i x_k}{\gamma_k f_k \wlr{x_k - x_i}}}^2},
\end{equation}
which is $\wfc{\Theta}{n^2 \epsilon}$ for $x_i = - \wfc{\cos}{i \pi / n}$, $k$ near $n$ and
$f_i = \wfc{\sin}{x_i}$.
We cannot use the central limit theorem to analyze the distribution of $s_{n,k}$,
because the variances of its terms do not satisfy Lindeberg's condition (see \cite{FELLER}, pg. 262). However,
it is reasonable to expect that $s_{n,k}$ will often assume values of the order of its standard deviation.
Thus, we can use \wref{var} as a rough estimate of the order of magnitude of the relative errors due to
rounding in the first formula for $t$ very close to the node $x_k$.
Therefore, Lemma \ref{lemChi} shows that
errors of order $n^2 \epsilon$, as listed in the second column of Tables \ref{tableNever} and \ref{tableNeverStd},
are to be expected.

{\small
\begin{table}[ht]
\caption{Values of $\wfc{s_{n,n-1}}{\hat{\wvec{x}}} / (n^2 \epsilon)$ for $\wfc{f}{t} = \wfc{\sin}{t}$.}
\centering
\begin{tabular}{c|r}
\hline\hline\\[-0.25cm]
 Nodes          & $\wfc{s_{n,n-1}}{\hat{\wvec{x}}} / (n^2 \epsilon)$ \\
\hline\\[-0.25cm]
$ 10^3 + 1$     & $0.039  $   \\
$ 10^4 + 1$     & $-0.044 $   \\
$ 10^5 + 1$     & $ 0.011 $   \\
$ 10^6 + 1$     & $ 0.013 $
\end{tabular}
\label{tableRKN}
\end{table}
}

The results concerning the stability of the second barycentric formula in section \ref{secSecond}
hold for general nodes. However, it is difficult to generalize Lemma \ref{lemChi} beyond the
Chebyshev points of the second kind. For other sets of nodes,
we must replace $\gamma_i / \gamma_k$ by the quotient $w_i/w_k$ appropriate for them
in expression \wref{snk}.
This complicates the analysis. For example, if the $x_i$ are equally spaced
then the ratio $w_{\wlr{n/2}} / w_{n-1}$ grows
exponentially with $n$ and the asymptotic arguments leading to Lemma \ref{lemChi} break down.
We emphasize therefore that the analysis in this section applies only to the
Chebyshev points of the second kind.
However, we hope our evidence illustrates
the  stability problems for interpolation of the first formula in general.
\subsection{The inconvenience of the first formula}
\label{subTheLastOne}
In the previous section, we saw that the first barycentric formula with simplified
weights and Chebyshev nodes of the second kind rounded as usual is unstable for large $n$,
  since we  may introduce relative errors of order $n^2 \epsilon$  by rounding the nodes.
We presented an experiment illustrating
this fact and a mathematical argument to explain it.
We conclude that we should not use simplified weights in combination with nodes rounded as usual for large $n$.
However, if we do decide to use better weights, then we face the problem of having to compute them first.
This is an expensive $\wfc{\Theta}{n^2}$ process.
 We don't believe improvements are possible by using Fourier techniques since we require accurate weights and our nodes are rounded.
As a result, there are two undesirable options:
we either compute the weights on the fly, at a significant $\wfc{\Theta}{n^2}$ cost,
or we precompute them and cope with the inconvenience of storing large tables
of weights, one for each $n$ we care about.

\section{Conclusions and future work}
\label{sec_the_end}
In summary, this article shows that the second barycentric formula as
considered by Salzer can be evaluated in a backward stable way.
For functions with derivatives of moderate size it also leads to small forward errors.
Furthermore, our article shows
that the first barycentric formula with nodes rounded as usual and simplified weights has stability problems.
Future work should explain how to implement the first formula with simplified weights in a stable way,
by handling the nodes differently,
and address issues regarding the stability of both formulae
that we did not consider here.

\section*{acknowledgements}
We thank M\'{a}rio Martinez for suggestions regarding the overall structure of this article,
Paulo Silva for helping with the numerical experiments, Andr\'{e} Camargo for a careful revision
of the technical details and Philippe Mettler for helping
with the English language. We thank both referees for their constructive
criticism and suggestions.

\appendix

\section{Proofs}
\label{apProofs}
This appendix contains proofs of our lemmas and theorems.

NOTE: in this appendix, we rely on our version of Stewart's notation \cite{HIGHAM_BOOK} to keep track of rounding errors.
Our notation is a slight generalization of Stewart's and we refer the reader to subsection \ref{subNotation} for details.

\subsection{Notation and conventions}
\label{subNotation}
Throughout the text we
used a hat to indicate the computed value of an expression, so that
$\hat{x}_i$ is the value we obtain by rounding the abstract $i$-th node
$x_i$.
The hat notation would however be cumbersome for large expressions.
Therefore, we follow \cite{HIGHAM_BOOK} and write
\[
\wrounde{a + b + c}
\hspace{2cm} \wrm{instead \ \ of} \ \hspace{2cm}
\widehat{a + b + c}.
\]
In other words, $\wrounde{\wrm{expression}}$ is the value we
obtain by evaluating $\wrm{expression}$ in floating
point arithmetic. Besides the hat, we use a version of Stewart's notation $\wst{\cdot}$.
This notation is based on the sets
\[
\werror{n} := \wset{ t \in \wrone{} \ | \
t = \prod_{i=1}^n \wlr{1 + \delta_i \epsilon}^{\sigma_i} \
\wrm{for \ some \ } \sigma_i \in \wset{-1,1} \ \wrm{and} \
-1 \leq \delta_i \leq 1}.
\]
As noticed by Stewart, when analyzing rounding errors it is convenient to
denote a generic element of $\werror{n}$ as $\wst{n}$ and say for example that
$\wrounde{x + y} = \wlr{x+y} \wst{1}$
to denote the assumption that there exists $z \in \werror{1}$ such that
$\wrounde{x + y} = \wlr{x+y} z$. We extend Stewart's notation in two ways.
First, as the reader can verify,
\begin{equation}
\label{xiRho}
\xi_\rho = \wset{t \in \wrone{} \ \ \wrm{with} \ \ \wlr{1 - \epsilon}^{\rho} \leq t \leq \wlr{1 - \epsilon}^{-\rho}}
\end{equation}
and we can use
this expression to define $\wst{\rho}$ for any $\rho \geq 1$, integer or not.
Second, sometimes it is necessary to give a name to the specific $\wst{\rho} \in \werror{\rho}$
we care about.
 We  use a subscript and denote this case with  $\wst{\rho}_a$.
 Throughout the proofs we follow this convention:
\begin{quote}
Every variable whose name is of the form $\wst{\rho}_a$ belongs
to the set $\werror{\rho}$. In particular, the equation
$a = b \ \wst{5}_c$
means not only that the left and right sides are equal but
also that $a/b \in \werror{5}$ when $b \neq 0$.
\end{quote}

Our floating point arithmetic is binary and there
 is no overflow or underflow. Therefore, if $m$ is an integer and $x$ is
 a floating point number then $2^m x$ is computed exactly. We often use that
\begin{equation}
\label{rhoToEps}
\rho \geq 1, \hspace{0.3cm}  \rho \epsilon \leq 0.01 \hspace{0.3cm} \wrm{and} \hspace{0.3cm}
 t \in \werror{\rho} \hspace{0.2cm}
\Rightarrow \hspace{0.2cm} \wabs{t - 1} \leq 1.01 \rho \epsilon.
\end{equation}
This can be proved by the same argument used to show the similar result on page
68 of \cite{HIGHAM_BOOK}.
The most convenient properties of Stewart's notation are:
\begin{equation}
\label{rhoProd}
\wst{\rho} \, \wst{\tau} \, = \, \wst{\rho \tau} \hspace{1cm} \wrm{and} \hspace{1cm} \frac{1}{\wst{\rho}_a} \, = \, \wst{\rho}_b,
\end{equation}
which follow directly from $\wref{xiRho}$. We use this version of the Standard Model on page 40 of \cite{HIGHAM_BOOK}:
\begin{equation}
\label{fModel}
\wrounde{x \ \wrm{op} \ y } = \wlr{x \ \wrm{op} \ y} \wst{1} \hspace{0.7cm} \wrm{for} \hspace{0.7cm}  \ \wrm{op} = +, -, * \ \ \wrm{and} \ \ /.
\end{equation}
We also assume that if $x \ \wrm{op} \ y$ is a floating point number then $\wrounde{x \ \wrm{op} \ y} = x \ \wrm{op} \ y $.

Finally, equation \wref{sumA} is a natural way to express the rounding errors in sums, which
does not require any background from the reader.
Therefore, using it instead of a criterion involving
Stewart's notation simplifies the presentation. However, there is a minor
incompatibility of $\wfc{O}{\epsilon}$ between these two notations. The next lemma shows that, to
leading order in $\epsilon$, we can say that $1 + \rho \delta \epsilon \in \werror{\rho}$
and presents a convenient set to which $1 + \rho \delta \epsilon$ belongs:
\begin{lem}
\label{lemFix} If $\rho \epsilon < 1 - 1/\sqrt{2}$ and $\wabs{\delta} \leq 1$ then $1 + \rho \delta \epsilon \in \werror{\rho \wlr{1 + \rho \epsilon}}$. \qed{}
\end{lem}

\subsection{An overview of the proofs of Theorem \ref{thmBackward} and Theorem \ref{thmMain}.}
\label{subOverview}
We now explain the ideas behind our theorems.
The first step is to write \wref{main} as $\wfc{b_n}{t;f} = \wfc{p}{t} / \wfc{q}{t}$ for
\begin{equation}
\label{pq}
\wfc{p}{t} := \sum_{i = 0}^n \frac{\gamma_i f_i}{t - x_i}
\hspace{1cm} \wrm{and} \hspace{1cm}
\wfc{q}{t} := \sum_{i =0}^n \frac{\gamma_i}{t - x_i}.
\end{equation}
It is possible to rewrite $\wfc{q}{t}$ as a sum of numbers with the same sign:
\begin{lem}
\label{lemDen}
If $n \geq 2$, $0 \leq k < n$, $x_0 = -1 < x_1 < \dots  < x_n = 1$,
 $t \in (x_k, x_{k+1})$ and the weights are as in \wref{theoGamma}
then we have four possibilities for $\wfc{q}{t}$, depending upon the parity of $k$ and $n$
\[
\begin{array}{r|cccc}
                              &     &      k \ \wrm{even} &   &  k \ \wrm{odd}  \\[0.2cm]
\hline\\
n = 2l \hspace{0.1cm} & \hspace{0.1cm} & \wfc{q}{t} =   \alpha   + \sum_{i = 1}^{l - 1} \xi_{2i    } + \psi   & \hspace{0.5cm} &
                                 \wfc{q}{t} =  - \wlr{ \beta    + \sum_{i = 1}^{l - 1} \xi_{2i + 1} + \omega } \\[0.2cm]
n = 2l + 1 \hspace{0.1cm} & \hspace{0.1cm} & \wfc{q}{t} =   \alpha   + \sum_{i = 1}^{l}     \xi_{2i    } + \omega & \hspace{0.5cm} &
                                  \wfc{q}{t} =  - \wlr{ \beta    + \sum_{i = 1}^{l-1}   \xi_{2i + 1} + \psi},
\end{array}
\]
where
\begin{equation}
\label{defAlpha}
\alpha  :=  \frac{1}{2 \wlr{1 + t}}, \hspace{0.6cm} \beta := \wfc{\eta}{t,x_1}, \hspace{0.6cm}
\psi    := \wfc{\eta}{-t, -x_{n-1}}, \hspace{0.6cm}
\omega  := \frac{1}{2 \wlr{1- t}},
\end{equation}
for
\begin{equation}
\label{eta}
\wfc{\eta}{y,z} :=  \frac{\wlr{2 + z} + y}{2 \wlr{y - z}\wlr{1 + y}} \hspace{1cm} \wrm{and} \hspace{1cm}
\xi_{i} :=  \frac{x_{i} - x_{i-1}}{\wlr{t - x_{i}} \wlr{t - x_{i-1}}}.
\end{equation}
Moreover, $\alpha$, $\beta$, $\psi$, $\omega$ and all the $\xi_i$ above are positive.
\end{lem}
Using Lemma \ref{lemFix} and Stewart's notation, we can show that
if all quantities $a_i$ are non negative then
$\hat{\sum}_{i = 0}^n a_i = \wst{1.01 \sigma_n} \sum_{i = 0}^n a_i$.
Therefore, we can evaluate this sum with high relative precision.
As a consequence, Lemma \ref{lemDen} shows that we can evaluate the
denominator of the second barycentric formula with high relative precision.
This is the key ingredient to prove Theorem \ref{thmBackward}.

The possibility of rewriting the denominator in order to avoid cancelations
was already noticed in Bos, De Marchi and Horman \cite{BOS}.
We have also found that Berrut \cite{BERRUT_POLE} and
Floater and Hormann \cite{FLOATER} have considered the question of absence
of poles in barycentric rational formulae, which is the mathematical essence of Lemma \ref{lemDen}.
Therefore, from a mathematical point of view,
we can say that Lemma \ref{lemDen} could be expected from the work of these authors.
However, we are not aware of references containing its explicit formula.
We are also unaware of any analysis of the numerical stability of the second
barycentric formula based on the work of these authors.

To get the bounds in Theorem \ref{thmMain} we
write the numerator $\wfc{p}{t}$ as a combination
with the positive coefficients $\xi_i$, $\alpha$, $\beta$, $\psi$ and $\omega$.
In other words, we write $\wfc{p}{t} = \sum \wfc{q_k}{t} \wfc{p_k}{t}$
so that the coefficients $q_k$ are positive, $q = \sum q_k$
and the $p_k$ can be evaluated accurately.
As a consequence, the second formula can be written
as a convex combination of the form
\begin{equation}
\label{convex}
\wfc{b_n}{t;f} = \sum \wlr{ \frac{\wfc{q_k}{t}}{\wfc{q}{t}}} \wfc{p_k}{t}.
\end{equation}
Comparing the expression for $\wfc{p}{t}$ in \wref{pq} and our
target $\wfc{p}{t} = \sum \wfc{q_k}{t} \wfc{p_k}{t}$ we derived the following lemma.
\begin{lem}
\label{lemNum}
If $n \geq 2$, $0 \leq k < n$, $x_0 = -1 < x_1 < \dots  < x_n = 1$,
 $t \in (x_k, x_{k+1})$ and the weights are as in \wref{theoGamma}
then we have four possibilities for $\wfc{p}{t}$, depending upon the parity of $k$ and $n$:
{\small
\[
\begin{array}{r|ccc}
                   & k \ \wrm{even} & k \ \wrm{odd}  \\[0.2cm]
\hline\\
n = 2l  & \hspace{0.1cm} \wfc{p}{t} =  \alpha f_0 + \sum_{i = 1}^{l-1} \xi_{2i} \phi_{2i} + \psi \psi_f
                              & \hspace{0.1cm} \wfc{p}{t} =  - \wlr{ \beta \beta_f + \sum_{i = 1}^{l - 1}  \xi_{2i + 1} \phi_{2i + 1} + \omega f_n}, \\[0.2cm]
n = 2l + 1  & \hspace{0.1cm} \wfc{p}{t} =  \alpha f_0 + \sum_{i = 1}^{l} \xi_{2i } \phi_{2i} + \omega f_n,
                              & \hspace{0.1cm} \wfc{p}{t} = - \wlr{ \beta \beta_f + \sum_{i = 1}^{l-1} \xi_{2i + 1} \phi_{2i + 1} + \psi \psi_f}\\
\end{array}
\]
}
where
\[
\beta_f := \wfc{\theta}{f_1, \ f_0, \ t, \ x_1} \hspace{1cm} \wrm{and} \hspace{1cm} \psi_f := \wfc{\theta}{f_{n-1}, \ f_n, \ - t, - \ x_{n-1}},
\]
for
\begin{equation}
\label{defTheta}
\wfc{\theta}{u, \ v, \ x, \ y} :=  \frac{2 \wlr{1 + x} u - \wlr{x - y} v}{\wlr{2 + y} + x},
\end{equation}
and
\[
\phi_i := \wfc{\varphi}{f_i, \ f_{i-1}, \ t, \ \frac{x_i + x_{i-1}}{2}, \ x_i, \ x_{i-1}}
\]
for
\begin{equation}
\label{defPhi}
\wfc{\varphi}{u, \ v, \ x, \ y, \ z, \ w} :=  \frac{u + v}{2} + \wlr{x - y} \frac{u - v}{z - w}.
\end{equation}
\wend{}
\end{lem}
 All the terms in Lemma \ref{lemNum}, except for the $\xi_i$,  are bounded
and can be evaluated in a backward and forward stable way.
It is not difficult to derive a proof of Theorem \ref{thmMain} from this lemma and \wref{convex}.

Finally, we note that
the denominator of $\xi$ in \wref{eta} can underflow
when $x_i = 0$, $n$ is large and $t$ is very close to $0$.
We ignore the issue here for two reasons: (i) it can be easily handled in practice with the introduction of \verb|if clauses| in the code
and (ii) it would make our analysis unnecessarily complicated.

\subsection{Proof of Theorem \ref{thmBackward}}
\label{subBackwardProof}
We assume that the numerator $\wfc{p}{t}$ in \wref{pq} is evaluated
naively. Lemma \ref{lemDen} shows that the denominator of the second barycentric
formula can be written as
\begin{equation}
\label{q_decomp}
q = s \wlr{q_{-} \, + \,  \sum_{k = 1}^{l_n} q_k \,  + \,  q_+},
\end{equation}
where $s \in \wset{-1,1}$ is a sign, $q_{-} \in \wset{\alpha, \beta}$, $q_{+} \in \wset{\psi, \omega}$
and $q_k = \xi_{i_k}$ for an appropriate index $i_k$ and
$l_n \in \wset{ \frac{n-3}{2}, \frac{n - 2}{2}, \frac{n-1}{2} }$.
The following lemma is about the rounded version
of the quantities $q$, $q_-$, $q_k$ and $q_+$ in equation \wref{q_decomp}.
\begin{lem}[Backward stability of the denominator]
\label{lemQHat}
If $n \geq 2$ and $2 + x_1$ and $2 - x_{n-1}$
are floating point numbers and $\sigma_n \epsilon \leq 0.01$ then,
for $m := \wfloor{ \wlr{n+1}/2 }$, the quantities in
\wref{q_decomp} can be evaluated so that
\[
\hat{q}_k = \hat{\xi}_{i_k} \ = \ q_k  \wst{5}_k, \hspace{0.4cm}
\hat{q}_- \ = q_- \wst{4}_-, \hspace{0.4cm} \hat{q}_+ \ =  q_+  \wst{4}_+
\]
\begin{equation}
\label{qMinusHat}
\hspace{0.3cm} \wrm{and} \hspace{0.4cm}
\hat{q} \ =  q \ \wst{5 + 1.01 \sigma_m}.
\end{equation}
\wend{}
\end{lem}

We prove this lemma in the next section. To prove Theorem \ref{thmBackward},
we note that there are no rounding errors in the multiplication $\gamma_i f_i$, because
$\gamma_i \in \wset{\pm 1/2, \pm 1}$. Therefore, using \wref{rhoProd} and \wref{fModel}
we can estimate the error in the
$i$-th term $p_i$ of the numerator as follows:
\[
\wrounde{\frac{\gamma_i f_i}{t - x_i}} = \frac{\gamma_i f_i}{t - x_i} \ \wst{2}_i.
\]
As a result, using \wref{sumA}, Lemma \ref{lemFix}  and \wref{rhoProd} we obtain
\[
\hat{p} = \sum_{i = 0}^n \frac{\gamma_i f_i}{t - x_i} \wst{2}_i \ \wlr{1 + \sigma_n \delta_i  \epsilon} =
\sum_{i = 0}^n \frac{\gamma_i f_i}{t - x_i} \wst{2 + 1.01 \sigma_n}_i.
\]
Lemma \ref{lemQHat},  \wref{rhoProd} and \wref{fModel} yield $m \leq (n + 1) / 2$ such that
\[
\wfc{\hat{b}_n}{t;f} = \wrounde{\frac{\hat{p}}{\hat{q}}} =
\frac{\hat{p}}{\hat{q}} \ \wst{1} =
\sum_{i = 0}^n \frac{\gamma_i f_i}{q \wlr{t - x_i}} \wst{8 + 1.01 \sigma_{m} + 1.01 \sigma_n}_i.
\]
We complete the proof of Theorem \ref{thmBackward} by taking
$\tilde{f}_i = f_i \, \wst{8 + 1.01 \sigma_{m} + 1.01 \sigma_n}_i$.
\qed{}

\subsection{Proof of Theorem \ref{thmMain}}
\label{subThmMain}
The previous proof already analyzed the
denominator $\wfc{q}{t}$, via Lemma \ref{lemQHat}.
Let us now look at the numerator $\wfc{p}{t}$
from the perspective of Lemma \ref{lemNum}.
 This lemma considers four cases but we can combine them
 as we did for $q$  in Lemma \ref{lemQHat}:
\begin{equation}
\label{pkMain}
\wfc{p}{t} = s \wlr{ q_- p_- \, + \,  \sum_{k = 1}^{l_n} q_k p_k \,  + \,  q_+ p_+},
\end{equation}
where $s$, $q_-$, $q_k$, $q_+$ and $l_n$ are defined just before
Lemma \ref{lemQHat}. The factors $p_k$ are defined as
\begin{equation}
\label{defPk}
p_k := \wfc{p_k}{f_{i_k}, f_{i_k-1}} := \wfc{\varphi}{f_{i_k}, \ f_{i_k-1}, \ t, \ \frac{x_{i_k} + x_{i_k-1}}{2}, \ x_{i_k}, \ x_{i_k   -1}},
\end{equation}
where $i_k$ is the same index used to define $q_k$ and the function $\varphi$ is
defined in \wref{defPhi}. Note our writing of
the $p_*$ both as constant and as a function of $f$, ignoring their dependency in $t$ and the $x_i$.
The term $p_-$ depends on the case we consider in Lemma \ref{lemNum}.
We may either have $p_- = f_0$ or
$p_- = \wfc{\theta}{f_1, \ f_0,\ t, \ x_1}$,
for $\theta$ defined in \wref{defTheta}.
Finally, we may either have $p_+ = f_n$ or
$p_+ = \wfc{\theta}{f_{n-1}, \ f_{n}, \ -t,\  -x_{n-1}}$.
The next two lemmas show that the $p_*$ in \wref{pkMain}
are bounded and backward and forward stable in the precise sense
explained in their conclusion:
\begin{lem}
\label{lemBetaF}
If formula \wref{defTheta} evaluated in the order indicated by the parentheses and
$2 + y$ is a floating point number then
$\wrounde{\wfc{\theta}{u,v,x,y,z}} = \wfc{\theta}{\tilde{u}, \tilde{v}, \ x , \ y, \ z}$, with
\begin{equation}
\label{backwardTheta}
\wabs{u - \tilde{u}} \leq 4.04 \wabs{u} \epsilon
\hspace{0.5cm} \wrm{and} \hspace{0.5cm}
\wabs{v - \tilde{v}} \leq 4.04 \wabs{v} \epsilon.
\end{equation}
We also have
\begin{equation}
\label{forwardTheta}
\wabs{\wrounde{\wfc{\theta}{u,v,x,y,z}} -\wfc{\theta}{u,v,x,y,z}} \leq  4.04 \times \wlr{2 \wabs{u} + \wabs{v}} \epsilon
\end{equation}
and if $-1 \leq y \leq x \leq 1$ then
\begin{equation}
\label{boundTheta}
\wabs{ \wfc{\theta}{u,v,x,y,z} } \leq 2 \wabs{u} + \wabs{v}.
\end{equation}
\wend{}
\end{lem}

\begin{lem}
\label{lemPhi}
Formula \wref{defPhi} evaluated in the order indicated by the parenthesis,
with $-1 \leq y,x \leq 1$, is backward stable with respect to $(u,v)$, in the sense that
\begin{equation}
\label{phiBackA}
\wrounde{\wfc{\varphi}{u,v,x,y,w, z}} =
\wfc{\varphi}{u + \wlr{a_{11} u + a_{12} v} \epsilon, \ v + \wlr{a_{21}  u + a_{22} v} \epsilon,\ x, \ y, \ w, \ z},
\end{equation}
with
\begin{equation}
\label{phiBackB}
\max \wset{ \wabs{a_{11}}, \wabs{a_{12}}, \wabs{a_{21}}, \wabs{a_{22}}} \leq 4.04.
\end{equation}
It is also forward stable, in the sense that
\begin{equation}
\label{forwardPhi}
\wabs{ \wrounde{\wfc{\varphi}{u,v,x,y,w, z}} - \wfc{\varphi}{u,v,x,y,w, z}} \leq 1.01
\wlr{\wabs{u} + \wabs{v} + 12 \wabs{\frac{u - v}{w - z}} } \epsilon.
\end{equation}
Moreover, $\varphi$ is bounded:
\begin{equation}
\label{boundPhi}
\wabs{\wfc{\varphi}{u,v,x,y,w,z}} \leq \frac{1}{2} \wlr{\wabs{u} + \wabs{v}} + 2 \wabs{ \frac{u - v}{z- w} }.
\end{equation}
\wend{}
\end{lem}

These lemmas, the stability of the denominator and the
equation \wref{qMinusHat} are the main ingredients to
obtain the backward error bound in equation \wref{backBoundMain}.
We prove \wref{backBoundMain} first and then continue to prove the forward error bound \wref{forwardBound}.

Equation \wref{pkMain} leads to
\[
\hat{p} = s \ \wrounde{  \wrounde{q_- p_-} \ + \  \sum_{k = 1}^{l_n} \wrounde{q_k p_k} \ + \ \wrounde{q_+ p_+}}
\]
and \wref{sumA} and Lemma \ref{lemFix} yield
\[
\hat{p} = s \wlr{ \wrounde{q_- p_-}  \wst{1.01 \sigma_m}_- \ + \ \sum_{k = 1}^{l_n} \wrounde{q_k p_k}  \wst{1.01 \sigma_m}_k \ + \
\wrounde{q_+ p_+}  \wst{1.01 \sigma_m}_+ },
\]
for $m := \wfloor{(n+1)/2} \geq l_n + 1$ (notice that the $\sigma_n$ in \wref{sumA} correspond to a sum of $n + 1$ terms.)
It follows from Lemma \ref{lemQHat} and equations \wref{rhoProd}  and \wref{fModel} that
\[
\wfc{\hat{b}_n}{t;f} = \wrounde{\frac{\hat{p}}{\hat{q}}} = \frac{\hat{p}}{\hat{q}}   \wst{1} \, = \,
\frac{\hat{p}}{q}   \wst{6 + 1.01 \sigma_m}.
\]
Therefore,
\[
\hfill
\wfc{\hat{b}_n}{t;f} = \frac{s}{q} \left( \wrounde{q_- p_-}  \wst{6 + 2.02 \sigma_m}_-  +
            \sum_{k = 1}^{l_n} \wrounde{q_k p_k} \wst{6 + 2.02 \sigma_m}_k  + \right.
\]

\begin{equation}
\label{qMain}
\Bigg.
\wrounde{q_+ p_+} \wst{6 + 2.02 \sigma_m}_+
\Bigg)
\hfill
\end{equation}
Equations \wref{qMinusHat}, \wref{rhoProd} and \wref{fModel} and the fact that $p_-$ is either $f_0$ or
$\wfc{\theta}{f_1, f_0, t, x_1}$ combined with Lemma \ref{lemBetaF} yield
\[
\wrounde{q_- p_-} = \hat{q}_- \hat{p}_- \wst{1} \, = \,
q_-  \wst{4} \wfc{p_-}{\tilde{f}_1,\tilde{f}_0} \wst{1} \, = \,
q_- \  \wfc{p_-}{\tilde{f}_1,\tilde{f}_0} \ \wst{5},
\]
for some $\tilde{f}$ with $\tilde{f}_0$ and $\tilde{f}_1$ such that
\begin{equation}
\label{fTilde0}
\wabs{\tilde{f}_0 - f_0} \leq 4.04 \wlr{\wabs{f_0} + \wabs{f_1}} \epsilon
\hspace{0.5cm} \wrm{and} \hspace{0.5cm}
\wabs{\tilde{f}_1 - f_1} \leq 4.04 \wlr{\wabs{f_0} + \wabs{f_1}} \epsilon.
\end{equation}
Analogously, we have
\[
\wrounde{q_+ p_+} = q_+ \wfc{p_+}{\tilde{f}_n,\tilde{f}_{n-1}} \wst{5},
\]
for $\tilde{f}$ with
\[
\wabs{\tilde{f}_{n-1} - f_{n-1}} \leq 4.04 \wlr{\wabs{f_{n-1}} + \wabs{f_n}} \epsilon
\hspace{0.4cm} \wrm{and} \hspace{0.5    cm}
\wabs{\tilde{f}_n - f_n} \leq 4.04 \wlr{\wabs{f_{n-1}} + \wabs{f_n}} \epsilon.
\]
Finally,
\[
\wrounde{q_k p_k} = \hat{q}_k \hat{p}_k \wst{1} \, = \,
q_k \wst{5} \wfc{p_k}{\tilde{f}_{i_k}, \tilde{f}_{i_k-1}}  \wst{1} \, = \,
\wfc{p_k}{\tilde{f}_{i_k}, \tilde{f}_{i_k-1}} \wst{6}
\]
for some $\tilde{f}$ with $\tilde{f}_{i_k - 1}$ and $\tilde{f}_{i_k}$ such that
\[
\wabs{\tilde{f}_{i_k - 1} - f_{i_k - 1}} \leq 4.04 \wlr{\wabs{f_{i_k - 1}} + \wabs{f_{i_k}}} \epsilon
\hspace{0.5cm} \wrm{and} \hspace{0.5cm}
\wabs{\tilde{f}_{i_k} - f_{i_k}}  \leq  4.04 \wlr{\wabs{f_{i_k - 1}} + \wabs{f_{i_k}}} \epsilon  .
\]
Combining the equations above with \wref{qMain} and using \wref{rhoProd} and \wref{fModel} we obtain
\[
\hfill
\wfc{\hat{b}_n}{t;f} =  \frac{s}{q} \Bigg( \
q_- \wfc{p_-}{\tilde{f}_1,\tilde{f}_0}  \wst{11 + 2.02 \sigma_m}_-
\]
\[
\left.
\ + \ \sum_{k = 1}^{l_n} q_k \wfc{p_k}{\tilde{f}_{i_k}, \tilde{f}_{i_k-1}}  \wst{12 + 2.02 \sigma_m}_k \ + \
q_+ \wfc{p_+}{\tilde{f}_{n-1}, \tilde{f}_{n}}  \wst{11 + 2.02 \sigma_m}_+
\right).
\hfill
\]
Since the $p_-$, $p_k$ and $p_+$ depend linearly in $f$ we can
collect the rounding errors in $\tilde{f}$, obtaining a vector $\hat{f}$
such that the second barycentric formula $b_n$ satisfies
\begin{equation}
\label{bMain}
\wfc{\hat{b}_n}{t;f} =  \frac{s}{q} \wlr{ q_- \wfc{p_-}{\hat{f}_1,\hat{f}_0} \ + \
\sum_{k = 1}^{l_n} q_k \wfc{p_k}{\hat{f}_{i_k}, \hat{f}_{i_k-1}} + q_+ \wfc{p_+}{\hat{f}_{n-1},\hat{f}_{n}} },
\end{equation}
as we now explain.
The step from $\tilde{f}$ to $\hat{f}$ is done as follows for $f_0$:
\[
\hat{f}_0 = \tilde{f}_0 \, \wst{11 + 2.02 \sigma_m}_-.
\]
The hypothesis on $\sigma_m$ implies
that $\wlr{11 + 2.02 \sigma_m} \epsilon \leq 0.01$.
Therefore, equation
\[
\hat{f}_0 - f_0 = \wlr{\tilde{f}_0 - f_0} \wst{11 + 2.02 \sigma_m}_- + f_0 \wlr{\wst{11 + 2.02 \sigma_m}_- - 1}
\]
and \wref{rhoToEps} and \wref{fModel} show that
\[
\wabs{\hat{f}_0 - f_0} \leq \wabs{\tilde{f}_0 - f_0} \times 1.01 + \wabs{f_0} \times 1.01 \times \wlr{11 + 2.02 \sigma_m}\epsilon.
\]
Combining this with bound \wref{fTilde0} we conclude that
\[
\wabs{\hat{f}_0 - f_0} \leq 1.01 \wlr{15.4 + 2.02 \sigma_m} \wlr{ \wabs{f_0} + \wabs{f_1} } \epsilon.
\]
The analysis of the remaining cases is analogous.
This finishes our analysis of the backward error.

Let us now analyze the forward error. We begin by noticing
that in all cases in which the function $\theta$ is evaluated in
Lemma \ref{lemNum} we have $1 < y \leq x \leq 1$. Therefore,
we can use Lemma \ref{lemBetaF} to bound $\varphi$ and then
bound $p_k$ in \wref{defPk}.
Let us start our analysis from equation \wref{qMain}, from which we obtain
\[
\wfc{\hat{b}_n}{t;f} =
\frac{s}{q}  \left( \hat{q}_- \hat{p}_-  \wst{7 + 2.02 \sigma_m}_-  \ + \
\sum_{k = 1}^{l_n} \hat{q}_k \hat{p}_k  \wst{7 + 2.02 \sigma_m}_k  \ + \ \right.
\]
\[
\Bigg. \hat{q}_+ \hat{p}_+ \wst{7 + 2.02 \sigma_m}_+ \Bigg).
\]
Using \wref{rhoProd}, \wref{fModel} and \wref{qMinusHat} we can rewrite this expression as
\[
\frac{s}{q}  \Bigg( q_- \hat{p}_-  \wst{11 + 2.02 \sigma_m}_-  \ + \Bigg.
\]
\begin{equation}
\label{almost}
\Bigg. \sum_{k = 1}^{l_n} q_k \hat{p}_k  \wst{12 + 2.02 \sigma_m}_k  \ + \
q_+ \hat{p}_+  \wst{11 + 2.02 \sigma_m}_+ \Bigg).
\end{equation}
The forward bounds \wref{forwardTheta} and \wref{forwardPhi} yield
\begin{eqnarray}
\label{forwardA}
\wabs{\hat{p}_- - p_-} & \leq & 12.12 \wnorm{f}_\infty \epsilon, \\
\nonumber
\wabs{\hat{p}_k - p_k} & \leq & 2.02 \wlr{ \wnorm{f}_\infty + 6 \wnorm{\Delta \! f}_\infty } \epsilon, \\
\nonumber
\wabs{\hat{p}_+ - p_+} & \leq & 12.12 \wnorm{f}_\infty \epsilon,
\end{eqnarray}
with the $\Delta \! f$ defined in the last line of the statement of Theorem \ref{thmMain}. Moreover,
\wref{boundTheta} and \wref{boundPhi} imply that
\begin{equation}
\label{lastIngredient}
\wabs{p_-} \leq 3 \wnorm{f}_\infty, \hspace{1cm}
\wabs{p_k} \leq  \wnorm{f}_\infty + 2 \wnorm{\Delta \! f}_\infty \hspace{1cm} \wrm{and} \hspace{1cm}
\wabs{p_+} \leq  3 \wnorm{f}_\infty.
\end{equation}
Equations \wref{rhoToEps} and \wref{lastIngredient} lead to these bounds:
\begin{eqnarray}
\label{dpMinus}
\wabs{\hat{p}_-  \wst{11 + 2.02 \sigma_m}_- - p_-} & \leq & 1.01 \wlr{45.12 + 6.06 \sigma_m} \wnorm{f}_\infty \epsilon,  \\
\nonumber
\wabs{\hat{p}_k  \wst{12 + 2.02 \sigma_m}_k - p_k} & \leq & 1.01 \left(15 \wnorm{f}_\infty + 37 \wnorm{\Delta  \! f}_\infty \right. \\
\label{dpk}
                          & & \left. + 2.02 \wnorm{f}_\infty \sigma_m + 4.04 \wnorm{\Delta  \! f}_\infty \sigma_m \right) \epsilon , \\
\label{dpPlus}
\wabs{\hat{p}_+  \wst{11 + 2.02 \sigma_m}_+ - p_+} & \leq & 1.01 \wlr{45.12 + 6.06 \sigma_m} \wnorm{f}_\infty \epsilon.
\end{eqnarray}
In fact, let us derive \wref{dpMinus}:
\[
\wabs{\hat{p}_-  \wst{11 + 2.02 \sigma_m}_* - p_-} \leq
\wabs{\hat{p}_- - p_- } \wst{11 + 2.02 \sigma_m}_* +
\]
\[
\wabs{p_-} \wabs{\wst{11 + 2.02 \sigma_m}_* - 1}
\]
\[
\leq 12.12  \wnorm{f}_\infty \epsilon \times 1.01  +
3 \wnorm{f}_\infty \times 1.01 \times \wlr{11 + 2.02 \sigma_m} \epsilon =
\]
\[
= 1.01 \wlr{45.12 + 6.06 \sigma_m} \wnorm{f}_\infty \epsilon.
\]
The derivation of the bound \wref{dpPlus} is analogous. Let us
then derive \wref{dpk}.
\[
\wabs{\hat{p}_k  \wst{12 + 2.02 \sigma_m}_k - p_k} \leq
\wabs{\hat{p}_k  - p_k } \wst{12 + 2.02 \sigma_m}_k
\]
\[
+ \wabs{p_k} \wabs{\wst{12 + 2.02 \sigma_m}_k - 1}
\]
\[
\leq 2.02 \wlr{ \wnorm{f}_\infty + 6 \wnorm{\Delta  \! f}_\infty } \epsilon \times 1.01 +
\wlr{\wnorm{f}_\infty + 2 \wnorm{\Delta  \! f}_\infty} \times 1.01 \times \wlr{12 + 2.02 \sigma_m} \epsilon.
\]
Using a calculator to handle the decimal numbers in this expression, one can
conclude that it is smaller than the right hand side of \wref{dpk}.
By combining the bounds in equations \wref{dpMinus}--\wref{dpPlus} we
conclude that the left hand side in each one of them is less than or equal to
\[
\zeta = 1.01 \wlr{45.12 \wnorm{f}_\infty + 37 \wnorm{\Delta  \! f}_\infty +
6.06 \wnorm{f}_{\infty} \sigma_m + 4.04 \wnorm{\Delta  \! f}_\infty \sigma_m} \epsilon.
\]
It follows from equation \wref{almost} and the positivity of the terms $q_*$ that
\[
\wabs{\wfc{\hat{b}_n}{t;f} - \wfc{b_n}{t;f}} \leq \wabs{\frac{1}{q}} \wlr{q_- + \sum_{k=1}^{l_n} q_k + q_+} \zeta.
\]
The sum in the numerator in this expression is equal to $\wabs{q}$ and we have proved the forward
bound $\wref{forwardBound}$.
\qed{}

\subsection{Proofs of the lemmas}
\label{subAuxiliary}
This section contains proofs of the lemmas up to this point.

%
%

{\bf Proof of Lemma \ref{lemSmallDf}.} This proof is left to the reader. \qed{}

{\bf Proof of Lemma \ref{lemChi}.}
We can write \wref{aTilde} as
\[
\wfc{\tilde{a}}{t,\hat{\wvec{x}}} = \wlr{-1}^n \frac{2^{n-1}}{n}
\wlr{ \prod_{i \neq k} \wlr{\frac{t - \hat{x}_i}{x_{k} - x_i}}}
\wlr{ \prod_{i \neq k} \wlr{x_{k} - x_i} } \wfc{s}{t,\hat{\wvec{x}}} \gamma_{k}  f_{k},
\]
with
\begin{equation}
\label{sdef}
\wfc{s}{t,\hat{\wvec{x}}} :=  1 + \wlr{t - \hat{x}_{k}} \sum_{i \neq k} \frac{\theta_{n,k.i}}{t - \hat{x}_i}
\hspace{0.7cm} \wrm{for} \hspace{0.4cm} \ \theta_{n,k,i} := \frac{\gamma_i f_i}{\gamma_k f_k}.
\end{equation}
Since we are using the simplified weights $\gamma_i$ in \wref{theoGamma} and the Chebyshev points of the second kind,
equation \wref{riez} yields
\[
\wlr{-1}^n \frac{2^{n-1}}{n} \wlr{ \prod_{i \neq k} \wlr{x_{k} - x_i} } \gamma_{k} = 1.
\]
It follows that
\begin{equation}
\label{tilA}
\wfc{\tilde{a}}{t,\hat{\wvec{x}}} =   f_k  \times \wfc{s}{t,\hat{\wvec{x}}} \times \wfc{g}{t} \times \wfc{h}{t,\hat{\wvec{x}}},
\end{equation}
where $g(t)$ is the $k$-th Lagrange fundamental polynomial
\[
\wfc{g}{t} := \prod_{i \neq k} \wlr{\frac{t - x_i}{x_k - x_i} }
\hspace{1cm} \wrm{and} \hspace{1cm}
\wfc{h}{t,\hat{\wvec{x}}} := \prod_{i \neq k} \wlr{1 - \frac{\hat{x}_i - x_i}{t - x_i}}.
\]
The bounds on $\alpha_{n,k}$ and $\delta_{n,k}$ in \wref{chiH} yield
\begin{equation}
\label{dxxk}
\wabs{t - x_k} \leq \frac{1}{24} \wabs{x_i - x_k}
\hspace{0.5cm} \wrm{and} \hspace{0.5cm}
\wabs{x_i - \hat{x}_i} \leq \frac{1}{24} \wabs{x_i - x_k}.
\end{equation}
By the triangle inequality,
$\wabs{t - x_i} \geq  \wabs{x_i - x_k} - \wabs{t - x_k} \geq 23 \wabs{x_i - x_k} / {24}$. Similarly,
$\wabs{t - \hat{x}_i} \geq \wabs{x_k - x_i} - \wabs{x_i - \hat{x}_i} - \wabs{t - x_k} \geq 11 \wabs{x_i - x_k} / {12}$.
These inequalities and \wref{dxxk} yield
\begin{equation}
\label{dxxi}
\wabs{\frac{x_i - x_k}{t - x_i}} \leq \frac{24}{23},
\hspace{1cm}
\wabs{ \frac{x_i - x_k}{t - \hat{x}_i}} \leq \frac{12}{11}
\hspace{1cm} \wrm{and} \hspace{1cm}
z := \wabs{\frac{\hat{x}_i - x_i}{t - x_i}} \leq \frac{1}{23}.
\end{equation}
The mean value theorem for $\wfc{f}{t} = \wfc{\ln}{1 + t}$ in the interval $[-z,z]$
yields $\wfc{\xi_i}{t,\hat{\wvec{x}}} \in [-z,z]$ such that
\[
\wfc{\ln}{1 - \frac{\hat{x}_i - x_i}{t - x_i}}
= - \frac{\hat{x}_i - x_i}{t - x_i}
\wlr{1 + \wfc{\mu_i}{t,\hat{\wvec{x}}} \frac{\hat{x}_i - x_i}{t - x_i}
}
\]
for
\begin{equation}
\label{muBound}
0 \leq \wfc{\mu_i}{t,\hat{\wvec{x}}} := - \frac{1}{2} \wdsf{f}{\wfc{\xi_i}{t,\hat{\wvec{x}} }} =
\frac{1}{2 \wlr{1 + \wfc{\xi_i}{t,\hat{\wvec{x}}}}^2} \leq
\frac{1}{2 \wlr{1 - 1/23}^2} = \frac{23^2}{2^3 \times 11^2}.
\end{equation}
This motivates the definition
\begin{equation}
\label{betanki}
\wfc{\beta_{n,k,i}}{t,\hat{\wvec{x}}} := \wlr{1 + \wfc{\mu_i}{t,\hat{\wvec{x}}} \frac{\hat{x}_i - x_i}{t - x_i}} \frac{x_k - x_i}{t - x_i}
\end{equation}
and leads to
\begin{equation}
\label{hk}
\wfc{h}{t,\hat{\wvec{x}}} = e^{- \sum_{i \neq k} \wfc{\beta_{n,k,i}}{t, \hat{\wvec{x}}} \frac{\hat{x}_i - x_i}{x_k - x_i}}.
\end{equation}
Moreover, the definitions in \wref{alphaDelta}, the hypothesis \wref{chiH}, \wref{dxxi} and the bound on $\mu$ above lead to
\[
\wabs{\wfc{\beta_{n,k,i}}{t,\hat{\wvec{x}}} - 1} \leq
\wabs{\frac{x_i - x_k}{t - x_i}} \wabs{\frac{x_k - t}{x_i - x_k}}  +
 \wfc{\mu_i}{t,\hat{\wvec{x}}} \wabs{\frac{x_i - x_k}{t - x_i}}^2 \wabs{\frac{\hat{x}_i - x_i}{x_i - x_k}} \leq
 \]
 \[
\leq \frac{24}{23} \wfc{\delta_{n,k}}{t} + \frac{72}{121} \wfc{\alpha_{n,k}}{\hat{\wvec{x}}}.
\]
Since $24/23 < 1.05$ and $72/121 < 0.6$, we have derived \wref{dBeta}.

We now write \wref{sdef} as
\begin{equation}
\label{sdefB}
\wfc{s}{t,\hat{\wvec{x}}} =
\wfc{s}{t,\wvec{x}} \wlr{1 - \frac{\wfc{p}{t,\hat{\wvec{x}}}}{\wfc{s}{t,\wvec{x}}}},
\end{equation}
for
\begin{equation}
\label{pk}
\wfc{p}{t,\hat{\wvec{x}}} :=
\sum_{i \neq k} \frac{\theta_{n,k,i} \wlr{\hat{x}_k - x_k}}{t - \hat{x}_i} -
\sum_{i \neq k} \frac{ \theta_{n,k,i} \wlr{t - x_k} \wlr{\hat{x}_i - x_i}}
{\wlr{t - \hat{x}_i}\wlr{t - x_i}}.
\end{equation}

The bound on $\wabs{f_k}$ in \wref{chiH} and $\wabs{\gamma_k} = 1$ yield
$\wabs{\theta_{n,k,i}} \leq 2$. The bound on $\delta_{n,k}$ in \wref{chiH}
and \wref{dxxi} lead to
\[
\wabs{\wfc{p}{t,\hat{\wvec{x}}}} \leq \frac{24}{11}
\wlr{
\sum_{i \neq k}\wabs{ \frac{\hat{x}_k - x_k}{x_i - x_k}} + \frac{24}{23} \times
\sum_{i \neq k} \wabs{\frac{t - x_k} {x_i - x_k}} \times
\sum_{i \neq k} \wabs{\frac{\hat{x}_i - x_i}{x_i - x_k}} }.
\]
Therefore, using the the definitions  in \wref{alphaDelta} and the bounds
in \wref{chiH} we obtain
\[
\wabs{\wfc{p}{t,\hat{\wvec{x}}}} \leq \frac{24}{11}
\wlr{\wfc{\delta_{n,k}}{\hat{x}_k} +
\frac{24}{23} \wfc{\delta_{n,k}}{t} \wfc{\alpha_{n,k}}{\hat{\wvec{x}}}}
\leq \frac{24}{11 \times 23}.
\]
The bound $\wabs{\theta_{n,k,i}} \leq 2$, that on $\delta_{n,k}$ in \wref{chiH}
and the first bound in \wref{dxxi} yield
\begin{equation}
\label{sBound}
\wabs{\wfc{s}{t,\wvec{x}} - 1} =
\wabs{
\sum_{i \neq k}  \theta_{n,k,i}
\frac{t - x_k}{x_i - x_k} \frac{x_i - x_k}{t - x_i} }
\leq 2 \times \wfc{\delta_{n,k}}{t} \times \frac{24}{23}
\leq 2 \times \frac{1}{24} \times \frac{24}{23} = \frac{2}{23}.
\end{equation}
Combining the last two bounds we obtain
\begin{equation}
\label{pOverS}
\wabs{ \frac{\wfc{p}{t,\hat{\wvec{x}}}}{\wfc{s}{t,\wvec{x}} }}
\leq \frac{23}{21} \times \frac{24}{11}
\wlr{\wfc{\delta_{n,k}}{\hat{x}_k} + \frac{24}{23} \wfc{\delta_{n,k}}{t}
\wfc{\alpha_{n,k}}{\hat{\wvec{x}}}
}
\leq \frac{24}{21 \times 11} < \frac{1}{9}.
\end{equation}
As before, the mean value theorem for $\wfc{f}{t} = \wfc{\ln}{1 + t}$ yields
$\wfc{\xi_i}{t,\hat{\wvec{x}}} \in [-1/9,1/9]$ such that
\begin{equation}
\label{lnp}
\wfc{\ln}{1 - \frac{\wfc{p}{t,\hat{\wvec{x}}}}{\wfc{s}{t,\wvec{x}} }}
= - \frac{\wfc{p}{t,\hat{\wvec{x}}}}{\wfc{s}{t,\wvec{x}} }
\wlr{
 1 +
\wfc{\tau}{t,\hat{\wvec{x}}}
\frac{\wfc{p}{t,\hat{\wvec{x}}}}{\wfc{s}{t,\wvec{x}} }
 }
\end{equation}
for
\begin{equation}
\label{tauk}
0 \leq \wfc{\tau}{t,\hat{\wvec{x}}} =
- \frac{1}{2} \wdsf{f}{\wfc{\xi_i}{t,\hat{\wvec{x}}}} =
\frac{1}{2 \wlr{1 + \wfc{\xi_i}{t,\hat{\wvec{x}}}}^2} \leq \frac{1}{2\wlr{1 - 1/9}^2}
= \frac{3^4}{2^7}
\leq 0.64.
\end{equation}
We then define
\begin{equation}
\label{defKappa}
\wfc{\kappa_{n,k,i}}{t,\hat{\wvec{x}}} :=
\frac{\theta_{n,k,i} \wlr{x_k - x_i}}{\wfc{s}{t,\hat{\wvec{x}}} \wlr{t - \hat{x}_i}}
\wlr{1 +
\wfc{\tau}{t,\hat{\wvec{x}}}
\frac{\wfc{p}{t,\hat{\wvec{x}}}}{\wfc{s}{t,\wvec{x}} }
 }
\end{equation}
and \wref{sdefB}--\wref{lnp} show that
\[
\wfc{s}{t,\hat{\wvec{x}}} =
\wfc{s}{t,\wvec{x}}
e^{\wfc{\ln}{1 -
\frac{\wfc{p}{t,\hat{\wvec{x}}}}{\wf{s}{t,\wvec{x}}}
}}
=
\wfc{s}{t,\wvec{x}}
e^{- \sum_{i \neq k} \wfc{\kappa_{n,k,i}}{t,\hat{\wvec{x}}} \frac{\hat{x}_k - x_k}{x_k - x_i}
+  \sum_{i \neq k} \wfc{\kappa_{n,k,i}}{t,\hat{\wvec{x}}} \frac{\wlr{t - x_k} \wlr{\hat{x}_i - x_i}}{\wlr{x_k - x_i} \wlr{t - x_i}}
}.
\]
Combining this equation with \wref{tilA} and \wref{hk} we obtain \wref{rnk}.

We now prove \wref{dKappa}. From $\wabs{\theta_{n,k,i}} \leq 2$,
\wref{pk}, \wref{pOverS} and \wref{tauk} we obtain
\[
\wabs{\wfc{\tau}{t, \hat{\wvec{x}}}
\frac{\wfc{p}{t, \hat{\wvec{x}}}}{\wfc{s}{t, \wvec{x}}}}
\leq \frac{3^4}{2^7} \times \frac{23}{21} \times \frac{24}{11}
\wlr{ \wfc{\delta_{n,k}}{\hat{x}_k} + \frac{24}{23}
\wfc{\delta_{n,k}}{t}\wfc{\alpha_{n,k}}{\hat{\wvec{x}}}}
\leq
\]
\[
\leq 1.6  \wlr{ \wfc{\delta_{n,k}}{\hat{x}_k} + \wfc{\delta_{n,k}}{t}\wfc{\alpha_{n,k}}{\hat{\wvec{x}}}}.
\]
From the definition of $\alpha_{nk}$ and $\delta_{n,k}$ in \wref{alphaDelta},
the hypothesis \wref{chiH}, $\wabs{\theta_{n,k,i}} \leq 2$, \wref{dxxi} and \wref{sBound} we get
\[
\wabs{\frac{x_k - x_i}
{\wfc{s}{t,\wvec{x}}\wlr{t - \hat{x}_i}} - 1} = \wabs{\frac{1}{\wfc{s}{t,\wvec{x}}}}
\wabs{
\frac{\hat{x}_i - x_i}{x_k - x_i}
\frac{x_k - x_i}{t - \hat{x}_i}
+
\frac{x_k - t}{x_k - x_i} \frac{x_k - x_i}{t - \hat{x}_i} +
\wlr{1 - \wfc{s}{t,\wvec{x}}}
} \hspace{3.5cm}
\]
\[
\leq
\frac{1}{21/23} \wlr{\wfc{\alpha_{n,k}}{\hat{\wvec{x}}} \frac{12}{11}
+ \wfc{\delta_{n,k}}{t} \frac{12}{11} + \wfc{\delta_{n,k}}{t} \frac{48}{23}
} \leq
\frac{12 \times 23}{11 \times 21} \wlr{\wfc{\alpha_{n,k}}{\hat{\wvec{x}}} + 3 \wfc{\delta_{n,k}}{t}}
\leq
\]
\[
\leq 1.2 \wlr{\wfc{\alpha_{n,k}}{\hat{\wvec{x}}} + 3 \wfc{\delta_{n,k}}{t}}.
\]
Using the inequality $\wabs{u v - 1} \leq \wabs{v} \wabs{u - 1} + \wabs{v - 1}$,
with $u = 1 + \tau p / s$ and $v = (x_k- x_i) / {s \wlr{t - \hat{x}_i}}$,
\wref{defKappa} and the last two bounds we obtain, when $\theta_{n,k,i} \neq 0$,
\[
\wabs{\wfc{\kappa_{n,k,i}}{t,\hat{\wvec{x}}} - \theta_{n,k,i}} \leq
2 \wabs{\frac{\wfc{\kappa_{n,k,i}}{t,\hat{\wvec{x}}}}{\theta_{n,k,i}} - 1} \leq
\]
\[
\leq 
4 \wlr{\wfc{\delta_{n,k}}{\hat{x}_k} + \wfc{\alpha_{n,k}}{\hat{\wvec{x}}} \wfc{\delta_{n,k}}{t} }
+ 2.4 \wlr{\wfc{\alpha_{n,k}}{\hat{\wvec{x}}} + 3 \wfc{\delta_{n,k}}{t}}
\]
and we have verified bound \wref{dKappa} for $\theta_{n,k,i} \neq 0$.
Since it clearly holds also for $\theta_{n,k,i} = 0$, we are done.
\qed{}

{\bf Proof of Lemma \ref{lemFix}.}
For $0 \leq t \leq 1 - 1/\sqrt{2}$, there exists $0 \leq \xi_t \leq 1$ such that
$\wfc{\ln}{1 - t} = -t - \xi_t t^2$. Thus,
\[
\frac{\wfc{\ln}{1 - \rho \epsilon}}{\wfc{\ln}{1 - \epsilon} } \leq \wlr{1 + \rho \epsilon} \rho.
\]
The next inequalities and the intermediate value theorem yield Lemma \ref{lemFix}:
\[
\wlr{1 - \epsilon}^{\rho\wlr{1 + \rho \epsilon}} \leq \wlr{1 - \epsilon}^{\frac{\wfc{\ln}{1 - \rho \epsilon}}{\wfc{\ln}{1 - \epsilon} }}
= 1 - \rho\epsilon \leq 1 + \rho \delta \epsilon \leq 1 + \rho \epsilon <
\]
\[
< \wlr{1 + \epsilon}^{\rho \wlr{1 + \rho \epsilon}} <
\wlr{1 - \epsilon}^{- \rho \wlr{1 + \rho \epsilon}}.
\]
\qed{}

{\bf Proof of Lemma \ref{lemDen}.} Let us start by showing that all quantities
that Lemma \ref{lemDen} claims to be positive are indeed positive.
Note
that $\xi_{2i} > 0$ when $k$ is even and
$\xi_{2i+1} > 0$ when $k$ is odd, because $\xi_j > 0$ if $t \not \in (x_{j-1}, x_{j})$ and
the way $k$ was chosen guarantees that $t \not \in [x_{2i - 1},x_{2i}]$ when $k$ is even
and $t \not \in [x_{2i},x_{2i + 1}]$ when $k$ is odd.

The parameters $\alpha$ and $\omega$ are clearly positive for $-1 < t < 1$.
The parameter $\beta$ appears only in the second column of the table in Lemma \ref{lemDen}.
Therefore, we only evaluate $\beta$ for $t \in (x_k,x_{k+1})$ with an odd $k$. This implies
that we only evaluate $\beta$ for $t > t_1$. By looking at the expressions for $\beta$
and $\eta$ in equations \wref{defAlpha} and \wref{eta}, we conclude that $\beta > 0$
for every $t$ that requires the use of $\beta$. Similarly, the parameter $\psi$ appears
only in the diagonal of the table in Lemma \ref{lemDen}. Therefore, we only need to
evaluate $\psi$ for $t \in (x_k,x_{k+1})$ for $k$ with the same parity as $n$.
This implies that $k + 1 \leq n - 1$ and $t < x_{n-1}$, or equivalently,
$-t > -x_{n-1}$.  By looking at the expressions for $\psi$
and $\eta$ in equations \wref{defAlpha} and \wref{eta}, we conclude that
$\psi > 0$ for every $t$ that requires the use $\psi$.
In summary, we only need to evaluate the parameters $\alpha$, $\beta$, $\psi$ and
$\omega$ in circumstances in which the resulting value is positive.

The verification of the algebraic identities in Lemma \ref{lemDen} is a tedious, error prone,
exercise
and is best evaluated with numerical code.
We did that and leave the corroboration
of our findings to the reader. \qed{}

{\bf Proof of Lemma \ref{lemNum}.} See the argument in the proof of Lemma \ref{lemDen}. \qed{}

%
%
{\bf Proof of Lemma \ref{lemQHat}.}
As an exercise in Stewart's notation, using \wref{rhoProd} and \wref{fModel} the reader
can derive
\begin{equation}
\label{xiHat}
\hat{\xi}_i = \xi_i  \wst{5}_i
\end{equation}
from \wref{eta}. Since $2 + x_1$ and $2 - x_{n-1}$ are computed exactly,
\wref{rhoProd}, \wref{fModel}
and \wref{defAlpha} imply that
\[
\hat{\alpha} = \alpha  \wst{2}, \hspace{0.5cm}
\hat{\beta}  = \beta  \wst{4}, \hspace{0.5cm}
\hat{\psi}   = \psi  \wst{4} \hspace{0.5cm} \wrm{and} \hspace{0.5cm}
\hat{\omega} = \omega \wst{2}.
\]
Therefore, $\hat{q}_- = q_-  \wst{4}$,  $\hat{q}_+  = q_+ \wst{4}$ and $\hat{q}_k = q_k \wst{5}$.
Equations \wref{sumA} and \wref{xiHat} with $q_* \geq 0$ imply that
\[
\hat{q} = s \, \wrounde{ \hat{q}_{-} \, + \,  \sum_{k = 1}^{l_n} \hat{q}_k \,  + \,  \hat{q}_+ }
=
 s \, \wlr{ \hat{q}_{-} \, + \,  \sum_{k = 1}^{l_n} \hat{q}_k \,  + \,  \hat{q}_+ } \wlr{1 + \sigma_{m} \delta \epsilon} =
q \wst{5} \wlr{1 + \sigma_{m} \delta \epsilon},
\]
for some $\delta \in [-1,1]$ and $m := \wfloor{\wlr{n+1}/2} \geq l_n + 1$ (notice that $\sigma_n$ in \wref{sumA} corresponds to a sum of $n + 1$ terms).
Lemma \ref{lemFix} yields $1 + \sigma_{m} \delta \epsilon = \wst{1.01 \sigma_{m}}$
and Lemma \ref{lemQHat} follows. \qed{}
%
%
{\bf Proof of Lemma \ref{lemBetaF}.}
Let us start by proving the bound  \wref{boundTheta}.
The derivative of $\theta$ with respect to $x$ is
\[
\theta_x = \frac{ \wlr{2 u - v} \wlr{2 + y + x} - \wlr{2 (1 + x) u - \wlr{x - y} v}}{\wlr{2 + y + x}^2}.
\]
The numerator of $\theta_x$ simplifies to $2 \wlr{u - v} \wlr{1 + y}$.
Therefore, it is independent of $x$. It follows that $\theta_x$ is always
zero or different from zero for all $x$. In either case, this implies
that the maximum absolute value of $\theta$ is achieved with $x = 1$ or $x = y$,
because $y \leq x \leq 1$ by hypothesis.
Evaluating $\theta$ in these extreme cases we get
\begin{equation}
\label{firstThetaBound}
\wfc{\theta}{u, v, y, y}  = \frac{2 \wlr{1 + y} u}{2 \wlr{1 + y} } = u \hspace{0.7cm} \wrm{and} \hspace{0.7cm}
\wfc{\theta}{u, v, 1, y}  = \frac{4 u - (1 - y) v}{3 + y}.
\end{equation}
In the second case, since $-1 \leq y \leq 1$,
\[
\wabs{\wfc{\theta}{u, v, 1, y}} \leq \frac{4 \wabs{u} + 2 \wabs{v} }{2} \leq 2 \wabs{u} + \wabs{v}.
\]
Combining this with \wref{firstThetaBound} we get bound \wref{boundTheta}.
Let us now verify the backward error bound \wref{backwardTheta}.
We have
\[
a := \wrounde{2 \wlr{1 + x} u} = 2 \wlr{1 + x} u  \wst{2}_a
\hspace{0.3cm} \wrm{and} \hspace{0.3cm}
b := \wrounde{ \wlr{x - y} v} =  \wlr{x - y} v  \wst{2}_b.
\]
Since  $2 + y$ is evaluated exactly, \wref{rhoProd} and \wref{fModel} yield
\[
c := \wrounde{\wlr{2 + y} + x} = \wlr{2 + y + x} \wst{1}_c.
\]
If follows that
$
d := \wrounde{a + b} = \wlr{2 \wlr{1 + x} u  \wst{2}_a + \wlr{x - y} v  \wst{2}_b} \wst{1}
$
and
\[
\wrounde{\wfc{\theta}{u,v,x,y}} = \wrounde{\frac{d}{c}} = \frac{d}{c} \wst{1}_\theta.
\]
We can rewrite the last equation as
\begin{equation}
\label{roundedTheta}
\hat{\theta} = \frac{1}{2 + y + x} \wlr{2 \wlr{1 + x} u  \wst{2}_a + \wlr{x - y} v  \wst{2}_b} \wst{2}_\theta.
\end{equation}
The backward error bound \wref{backwardTheta} follows from
combining \wref{rhoToEps} with the following version of equation \wref{roundedTheta},
\[
\hat{\theta} = \frac{1}{2 + y + x} \wlr{2 \wlr{1 + x} \tilde{u}  + \wlr{x - y} \tilde{v} },
\]
where $\tilde{u} = u  \wst{2}_a \wst{2}_\theta$ and $\tilde{v} =  v  \wst{2}_b \wst{2}_\theta$.

Finally, let us verify the forward error bound \wref{forwardTheta}. The argument involving
the derivative $\hat{\theta}$ used to bound the maximum of $\hat{\theta}$
shows that the forward error $\wabs{\hat{\theta} - \theta}$
will not decrease if we replace $x$ by $y$ or $1$.
Therefore, we only need to consider
these two cases, which,  according
to equation \wref{roundedTheta},  lead to
\begin{eqnarray}
\nonumber
\theta - \hat{\theta} & = &  u  \wlr{1 - \wst{2}_a \wst{2}_\theta}. \\
\nonumber
\theta - \hat{\theta} & = &  \frac{1}{3 + y} \wlr{4 u + \wlr{1 - y} v} -
\frac{1}{3 + y} \wlr{4 u  \wst{2}_a + \wlr{1 - y} v  \wst{2}_b} \wst{2}_\theta \\
\nonumber
& = &  \frac{1}{3 + y} \wlr{4 u \wlr{1 -   \wst{2}_a \wst{2}_\theta} + \wlr{1 - y} v \wlr{1 - \wst{2}_b \wst{2}_\theta}}.
\end{eqnarray}
Since $-1 < y < 1$, in the second case the right hand side of
\[
\wabs{\theta - \hat{\theta}} \leq
\frac{1}{\wabs{3 + y}} \wlr{4 \wabs{u} \wabs{1 -   \wst{2}_a \wst{2}_\theta} + \wabs{1 - y} \wabs{v} \wabs{1 - \wst{2}_b \wst{2}_\theta}}.
\]
is maximized by taking $y = -1$ and we have
\[
\wabs{\theta - \hat{\theta}} \leq 2 \wabs{ u } \wabs{1 -   \wst{2}_a \wst{2}_\theta}
+ \wabs{v} \wabs{1 -  \wst{2}_b \wst{2}_\theta}.
\]
The forward bound \wref{forwardTheta} follows from \wref{rhoToEps}.
\qed{}
%
%
{\bf Proof of Lemma \ref{lemPhi}.} Equations \wref{rhoProd} and \wref{fModel} yield
\[
a := \wrounde{\frac{u + v}{2}} = \frac{u + v}{2} \wst{1}
\hspace{0.6 cm} \wrm{and} \hspace{0.6cm}
b := \wrounde{\wlr{x - y} \wlr{\frac{u - v}{w - z}}} =  
\]
\[
= \wlr{x - y} \wlr{\frac{u - v}{w - z}}  \wst{5}.
\]
Therefore,
\begin{equation}
\label{backPhi}
\wrounde{\wfc{\varphi}{u,v,x,y,w, z}} = \wrounde{a + b} = \wlr{a + b}  \wst{1}
\; = \; \frac{u + v}{2} \wst{2}_a + \frac{x - y}{w - z} \wlr{u - v} \wst{6}_b.
\end{equation}
We want to find $\tilde{u} = u + \wlr{a_{11} u + a_{12} v} \epsilon$ and
$\tilde{v} = v + \wlr{a_{21} u + a_{22} v} \epsilon$
such that \wref{backPhi} is equal to
\[
\frac{1}{2} \wlr{\tilde{u} + \tilde{v}} + \frac{x - y}{w - z}  \wlr{\tilde{u} - \tilde{v}}.
\]
We can achieve this goal by solving the following linear system for the variables $\tilde{u}$, $\tilde{v}$:
\[
\tilde{u} + \tilde{v}  =  \wlr{u + v} \wst{2}_a \hspace{0.5cm} \wrm{and} \hspace{0.5cm}
\tilde{u} - \tilde{v}  =  \wlr{u - v} \wst{6}_b.
\]
Its solution is
\begin{eqnarray}
\nonumber
\tilde{u} & = & u + \wlr{\frac{\wstx{2}_a + \wstx{6}_b}{2} - 1} u + \wlr{\frac{\wstx{2}_a - \wstx{6}_b}{2} } v \\
\nonumber
\tilde{v} & = & v + \wlr{\frac{\wstx{2}_a - \wstx{6}_b}{2}} u + \wlr{ \frac{\wstx{2}_a + \wstx{6}_b}{2} - 1} v.
\end{eqnarray}
It follows that
\[
\begin{array}{ll}
a_{11}  =  \frac{1}{\epsilon} \wlr{\frac{\wstxx{2}_a + \wstxx{6}_b}{2} - 1} & \hspace{1cm} a_{12}  =  \frac{1}{\epsilon} \wlr{\frac{\wstxx{2}_a - \wstxx{6}_b}{2}}, \\[0.1cm]
a_{21}  =  \frac{1}{\epsilon} \wlr{\frac{\wstxx{2}_a - \wstxx{6}_b}{2}}     & \hspace{1cm} a_{22}  =  \frac{1}{\epsilon} \wlr{\frac{\wstxx{2}_a + \wstxx{6}_b}{2} - 1}.
\end{array}
\]
Equation \wref{rhoToEps} shows that
the $a_{ij}$ above satisfy $\wabs{a_{ij}} \leq 4.04$ and we have proved
\wref{phiBackA}--\wref{phiBackB}.

The definition \wref{defPhi} and the triangle inequality yield \wref{boundPhi}:
\[
\wabs{\wfc{\varphi}{u,v,x,y,z,z}} \leq \frac{1}{2} \wlr{\wabs{u} + \wabs{v}}
+ \wabs{x - y} \wabs{\frac{u - u}{z - w}}.
\]
Now, $\wabs{x - y} \leq 2$ because $-1 \leq x, y \leq 1$ and \wref{boundPhi} follows.

Finally, to verify the forward error bound \wref{forwardPhi}, notice that
equation \wref{backPhi} yields
\[
\hat{\varphi}  - \varphi =
\frac{u + v}{2} \ \wlr{\wst{2}_a - 1} + \frac{x - y}{w - z} \wlr{u - v} \ \wlr{\wst{6}_b - 1}.
\]
Equation \wref{rhoToEps} and the bound $\wabs{x - y} \leq 2$ yield
\[
\wabs{\hat{\varphi}  - \varphi} \leq
\frac{\wabs{u} + \wabs{v}}{2} \ 2.02 \epsilon + 2 \wabs{\frac{u - v}{w - z}} 6.06 \epsilon.
\]
This is equivalent to bound \wref{forwardPhi} and we are done.
\qed{} 
\section{Experimental settings}
\label{apExperiments}
All experiments were performed on a Intel core i7
processor running Ubuntu 12.04. The code was written in \verb C++  and
compiled with g++4.7.0, with options
\verb -std=c++11  \verb -O3  \verb -march=corei7  \verb -mavx  \verb -Wall \verb -pthread .
Ultimate performance was not our concern, neither with
our methods nor with the others.

Errors were measured
with quadruple precision
 (1 bit sign, 128 bit unsigned mantissa and 32 bit signed
exponent).
 By $\wabs{\hat{y} - y}$ we mean the error in a computed value $\hat{y}$ corresponding to an exact value $y$.
 By relative error
we mean $\wabs{\hat{y} - y}/\wabs{y}$. When $y =0$ and
$\hat{y} \neq 0$ we say that the relative error is $\infty$, as in some tables in section 3.1.
We computed several results with even higher precision and they agreed with the quadruple precision ones.
The multi precision computations were performed with the MPFR library
\cite{MPFR}.
 Processing was timed
using the cpu time clock from the boost library
\cite{BOOST_SITE,BOOST_DEMING}, which is represented by the class
\verb boost::chrono::process_user_cpu_clock .
This is an accurate timer and it considers only the time taken
by the process one is concerned with.

The first formula was implemented using the products computed in appendix \ref{apOverflow}
and
 the sums were computed in the natural way.
The Chebyshev points of the second kind
were evaluated using the expression $x_i = \wfc{\sin}{\wlr{\wlr{2 i - n} \pi}/{2n}}$,
which is mathematically equivalent to the usual formula
$- \wfc{\cos}{i \pi / n}$ but has better numerical properties.
The $f_i$ were obtained by evaluating $\wfc{f}{\hat{x}_i}$  in quadruple precision
and then rounding the result to \verb double .
In tables \ref{tableNever} and \ref{tableNeverStd}, the first formula with simplified weights uses the
values in \wref{theoGamma} as the weights and the weights in \wref{firstNu} were obtained by
evaluating \wref{firstNu} in double precision. The Stable interpolant implements the formula in lemmas
\ref{lemDen} and \ref{lemNum}. The rounded nodes $\hat{x}_i$ were computed
using the code in appendix \ref{apRoundingNodes}. The $f_i$ were obtained by evaluating
$\wfc{f}{\hat{x}_i}$ in quadruple precision and then rounding the result to \verb double . Sums were
computed naively.

\section{Rounding the Chebyshev points of the second kind}
\label{apRoundingNodes}
This appendix  explains how to round the Chebyshev points of the second kind,
$x_i = - \wfc{\cos}{i \pi / n}$, to obtain $\hat{x}_i$ as
required by Theorem \ref{thmBackward} and Theorem \ref{thmMain}.
We present the \verb C++  code at the end of this appendix.
The code below is real, but to understand it
it is better to think in terms the idealized \verb double  and
\verb long  \verb double  numbers of the form
\[
D = \wset{x = 2^{\beta} m, \ \wrm{for \ integers } \ \beta, m \ \ \wrm{with} \ \beta \geq -1000
\ \wrm{and} \ 2^{52} \leq m < 2^{53} },
\]
\[
L = \wset{w = 2^{\beta} m, \ \wrm{for \ integers } \ \beta, m \ \ \wrm{with} \ \beta \geq -1000
\ \wrm{and} \ 2^{63} \leq m < 2^{64} }.
\]
This idealization takes  underflow into account by considering exponents
greater than $-1000$.
The smallest positive node is at least
$\wfc{\sin}{\pi/2n}$ and this allows us to handle $n \leq 10^9$.
We can safely ignore overflow since the nodes are small.
We only consider  positive nodes, because the others are $0$ or symmetrical.

Regarding the numbers $x, y \in D$ and $z,w \in L$, the code assumes the following:\\[-0.25cm]
\begin{enumerate}
\item Sums are rounded to nearest, i.e. $\widehat{x + y}$
is the number closest to $x + y$. Ties are broken arbitrarily.\\[-0.2cm]
\item If $\hat{x}_i \in L$ is the rounded version of $x_i$ then
$\wabs{\hat{x}_i - x_i} \leq 0.53 \times 2^{-52} \, \min \wset{\wabs{x_i},\wabs{\hat{x}_i}}$.
\item The function \verb sinl  is monotone in $[0,\pi/2] \cap L$, i.e., if $a < b \in L$ then
\verb sinl(a)  \verb <=  \verb sinl(b)  .
\item If $z = 2^\beta m$ then \verb frexpl(z,&exp)  returns $2^{-64} m$ and sets
\verb exp $ = 64 + \beta$.\\[-0.2cm]

\item if $x = 2^\beta m$ with $\beta > -900$ and
$k$ is an integer with $k > -100$ then \verb scalbln(x,k)  returns $2^{k} x$.\\[-0.2cm]

\item If $z = 2^\beta m$ and $k$ is a positive integer then \verb scalblnl(y,k)  returns $2^{k} z$.\\[-0.2cm]
\item If $i$ is a 64 bit unsigned integer such that $0 \leq i \leq 2^{53}$
then \\ \verb static_cast<double>(i)  returns an element of $D$ which is
mathematically equal to $i$.\\[-0.2cm]
\item If $z = 2^{0} m$ then \verb static_cast<uint64_t>(z)  returns $m$.\\[-0.2cm]
\end{enumerate}
Except for item 2, we believe
these assumptions are satisfied by most modern \verb C++  compilers
for most processors in use today. The Visual C++ compiler is a notorious exception with regard to item 2, because
it does not support the type \verb long  \verb double .
In item 2 we ask for relative errors smaller than $0.53 \times 2^{-52}$ while computing
the nodes using \verb long  \verb double  arithmetic. Since this arithmetic
has machine epsilon equal to $2^{-64}$, one may think that the computed
nodes would have relative errors much smaller than $0.53 \times 2^{-52}$.
In fact, our experiments with gcc indicate a maximum relative error of
$0.5003 \times 2^{-52}$ in the \verb long  \verb double  nodes.
 Our hypothesis still holds but only by a narrow margin.
The size of these relative errors illustrate that the evaluation of
the nodes involve more than a simple call to the $\sin$ function.
We must take into account that $\pi$ is also rounded.

The first step to build the rounded nodes $\hat{x}_i$ is to take
$y_n = 1$ and, for $n/2 < i < n$,
take $y_i \in L$ such that
$\wabs{x_i - y_i} \leq 0.53 \min\wset{\wabs{y_i},\wabs{x_i}}$.
Next we take $\hat{x}_n = 1$ and for $n/2 < i < n$ we choose $\hat{x}_i$ as follows:
\begin{itemize}
 \item If $y_i$ and $y_{i+1}$ have the same exponent then we choose $\hat{x}_i$ as the element
 of $D$ with an even mantissa closest to $y_i$.
 \item Otherwise, $\hat{x}_i$ is chosen as the element of $D$ with a mantissa multiple of
 four closest to $y_i$.
\end{itemize}
If $i = n/2$ then we take $\hat{x}_i = 0$ and if $0 \leq i < n/2$ then we take
$\hat{x}_i = - \hat{x}_{n - i}$. The correctness of this procedure
is assured by the following lemma, which we proved at the end of this appendix.

\begin{lem}
\label{lemRound}
For $n \leq 10^9$, $x_i = - \wfc{\cos}{i \pi / n}$ and $y_i \in L$ such that
\[
\wabs{x_i - y_i} \leq 0.53 \min \wset{\wabs{x_i},\wabs{y_i}} 2^{-52},
\] the $\hat{x}_i$
in the previous paragraph satisfy
$\wabs{\hat{x}_i - x_i} \leq 2.54 \times 2^{-52} \wabs{x_i}$ and 
\[
\wset{2 + \hat{x}_1, 2 - \hat{x}_{n-1}, \hat{x}_i + \hat{x}_{i+1}, i = 0, \dots, n-1} \subset D.
\]
\end{lem}

Here is the \verb C++  code implementing the ideas described above:
{
\small
\begin{verbatim}
double roundToEven(long double x) { // x is positive
  int32_t exp;
  long double xr = frexpl(x, &exp);
  uint64_t ix = static_cast< uint64_t >( scalblnl(xr, 64));
  ix >>= 11;
  if( ix & 0x1 ) ++ix;
  double rx = static_cast<double>(ix);
  return scalbln(rx, exp + 11 - 64);
}

double roundToMultipleOfFour(long double x) { // x is positive
  int32_t exp;
  long double xr = frexpl(x, &exp);
  uint64_t ix = static_cast< uint64_t >( scalblnl(xr, 64));
  ix >>= 11;
  switch( ix & 0x3 ) {
    case 1: --ix; break;
    case 2: ix += 2; break;
    case 3: ++ix; break;
  }
  double rx = static_cast<double>(ix);
  return scalbln(rx, exp + 11 - 64);
}

void roundedChebyshevNodes(double* nodes, uint32_t degree) {
  double* end = nodes + degree;
  *end   = 1.0;
  *nodes = -1.0;
  if( degree < 2 ) {
    return;
  }

  double xn;
  long double power = 1.0l;
  long double piOver2N = M_PIl / (2 * degree);
  while( ++nodes < --end ) {
    long double sn = sinl( (degree -= 2) * piOver2N);
    if( sn < power ) {
      power = scalbln(power,-1);
      xn = roundToMultipleOfFour(sn);
    }
    else  {
      xn = roundToEven(sn);
    }
    *end = xn;
    *nodes = -xn;
  }
  if( nodes == end ) *nodes = 0.0;
}
\end{verbatim}
}

{\bf Proof of Lemma \ref{lemRound}}. We consider only positive $x_i$, or $i > n/2$.
The bound on $\wabs{\hat{x}_i - x_i} $ follows from the hypothesis
$
\wabs{x_i - y_i} \leq 0.53 \min  \wset{\wabs{x}_i,\wabs{y_i}}2^{-53}
$
and the inequality $\wabs{\hat{x}_i - y_i} \leq 2^{-51} \wabs{y_i}$,
which is a consequence of the way we round $y_i$ to obtain $\hat{x}_i$.
Since $5 \leq n \leq 2^9$ we have $1/2 < \wfc{\sin}{3\pi/5} \leq y_{n-1} \leq 1$.
If $y_{n - 1} = 1$ then $y_{n-1} = \hat{x}_{n-1} = 1$ and
$2 - \hat{x}_{n-1} = 1 \in D$. Otherwise, the mantissa of $\hat{x}_{n-1}$ is
rounded to a multiple of four and we have
$\hat{x}_{n-1} = 2^{-53}{4m}$ with $2^{48} \leq m < 2^{49}$. It follows
that $2 - \hat{x}_{n-1} = 2^{-52} \wlr{2^{53}  -2 m} \in D$.
Therefore, we proved that $2 - \hat{x}_{n - 1} \in D$.
We now show that $\hat{x_i} +\hat{x}_{i + 1} \in D$ for $n/2 < i < n$.
The mantissas of the $\hat{x}_i$ are even and $\hat{x}_i \leq \hat{x}_{i+1}$, thus
\[
\hat{x}_i = 2^\beta 2u
\hspace{0.2cm} \wrm{with} \hspace{0.2cm}
2^{51} \leq u < 2^{52}
\hspace{0.5cm}
\hat{x}_{i+1} = 2^{\beta + \delta} 2v
\hspace{0.2cm} \wrm{with} \hspace{0.2cm}
 \delta \geq 0
\hspace{0.2cm} \wrm{and} \hspace{0.2cm}
 2^{51} \leq v < 2^{52}.
\]
If $\delta = 0$ then $\hat{x}_i + \hat{x}_{i+1} = 2^{\beta+1}\wlr{u+v} \in D$,
because $2^{52} \leq u+v < 2^{53}$. Therefore, we can assume that $\delta > 0$.
If $\hat{x}_{i+1} = 2^{\beta} 2^{53}$ then
$\hat{x}_i + \hat{x}_{i+1} = 2^{\beta+1} \wlr{u+2^{52}} \in D$. As a result,
we only need to concern ourselves with
$\hat{x}_{i+1} \geq 2^\beta \wlr{2^{53} +2}$.
Thus, $y_{i+1} \geq 2^\beta \wlr{2^{53} + 1} \geq 2^{\beta - 10}
\wlr{2^{63} +2^{10}}$. On the other hand,
$y_i \leq 2^\beta \wlr{2^{53} - 1} = 2^{\beta - 11} \wlr{2^{64} - 2^{11}}$.
It follows that $y_{i+1}$ has exponent at least $\beta - 10$
and $y_i$ has exponent at most $\beta - 11$. Therefore, these exponents are
different and, by construction, we only need to consider $\hat{x}_i$ given by
$\hat{x}_i = 2^{\beta} 4w$ with $2^{50} \leq w < 2^{51}$.
We now show that
\begin{equation}
\label{eq69}
\hat{x}_{i+1} \leq 2^{\beta +2} \wlr{2^{52} +2^{51}}.
\end{equation}
As a first step we show that $x_{i+1} \leq 3 x_i$. In fact,
since $x_i \geq \wfc{\sin}{\pi /{2n}}$ for all $i > n/2$, we have
\[
x_{i+1} = \wfc{\sin}{\frac{\wlr{2 i - n} \pi}{2n} + \frac{\pi}{n}}
\leq \wfc{\sin}{\frac{\wlr{2 i - n} \pi}{2n}} +\wfc{\sin}{\frac{\pi}{n}} =
\]
\[
= x_i +2 \wfc{\cos}{\frac{\pi}{2n}} \wfc{\sin}{\frac{\pi}{2n}} \leq 3 x_i.
\]
Now, by the way $\hat{x}_i$ was rounded,
$y_i \leq 2^\beta \wlr{2^{53} - 2}$
and
\[
y_{i+1} \leq x_{i+1} \wlr{1+ 0.53\times 2^{-52}}
\leq 3 x_i \wlr{1+0.53 \times2^{-52}} \leq 3 y_i
\wlr{1+0.53 \times 2^{-52}}^2
\]
\[
\leq 3\times 2^\beta \wlr{2^{53} - 2}
\wlr{1+1.06 \times 2^{-52} + 2^{-105}}
\leq 3\times 2^\beta
\wlr{2^{53} + 0.13} = 
\]
\[
= 2^\beta \wlr{3 \times2^{53} +0.39}.
\]
Equation \wref{eq69} follows from the way $\hat{x}_{i+1}$ is built
from $y_{i+1}$. This leaves us with three possibilities:
\begin{itemize}
\item[(a)] $\delta = 1$ and $2^{52} \leq 2v < 2^{53}- 2 w$,
\item[(b)] $\delta = 1$ and $2^{53}- 2w \leq 2v < 2^{53}$,
\item[(c)] $\delta = 2$ and $2v \leq 2^{52} +2^{51}$.
\end{itemize}
In case (a) we have
$\hat{x}_i +\hat{x}_{i+1} = 2^{\beta+1}\wlr{2w+2v} \in D$ because
$2^{52} \leq 2v+2w < 2^{53}$.
In case (b),
$\hat{x}_i + \hat{x}_{i+1} = 2^{\beta+2} \wlr{v+w} \in D$,
 because $2^{52} \leq v+w < 2^{53}$ in this case.
Finally, in case (c), $\hat{x}_i+\hat{x}_{i+1} = 2^{\beta+2} \wlr{2v+w} \in D$
because in this case
$2^{52} \leq 2v+w < 2^{52} +2^{51} +2^{51} = 2^{53}$.
\qed{}

\section{A robust auxiliary C++ function to evaluate the
first barycentric formula}
\label{apOverflow}
This appendix presents a \verb C++  function to evaluate the product
\wref{firstProduct} without spurious overflow or underflow and
with no sacrifice in terms of performance or accuracy. The function can be used
 under the following conditions:
\begin{itemize}
\item $-1 \leq t \leq 1$ and $1 \leq n \leq 10^9$.\\[-0.25cm]
\item The \verb xi[i]   are close to the exact nodes
$x_i = - \wfc{\cos}{i \pi / n}$, so that $|x_i - $ \verb xi[i]  $| \leq 0.99 \wabs{x_i}$.
\end{itemize}

The function evaluates most products in groups of twenty.
Once it is done with a group, it uses \verb frexp  to scale the product by a power of two,
using the variable \verb exp  to keep track of any scaling.
The potentially small factors $t - x_i$ with $x_i$ very close to $t$
are handled individually with one \verb frexp  per factor.
It is possible to prove that,  by handling the six nearest neighbors in each side of $t$  individually,
we obtain a correct product for $n \leq 10^9$. We do not show the proof for the sake of
brevity, but the reader can write a simple routine to verify that
products of the form $0.99^{20} \wlr{x_i - x_{i+6}} \wlr{x_i - x_{i+7}} \dots \wlr{x_i - x_{i+26}}$
do not underflow for $26 \leq i + 26 \leq n \leq 10^9$. This is
what one needs to prove that the code below works.

{\small
\begin{verbatim}
double firstProduct(double t, double const* xi,  int32_t n)
{
  if( t <= -1 ) return (t < -1.0) ? NAN : 0; // We do not allow t < -1.
  if( t >= 1  ) return (t >  1.0) ? NAN : 0; // We do not allow t > 1.
  if( (n < 1) || (n > 1000000000) ) return NAN;
  // finding min such that xi[min] <= t < xi[min + 1]
  int32_t min = 0;
  int32_t max = n;
  while( min + 1 < max ) { // keep the invariant xi[min] <= t < xi[max]
    int32_t middle = (min + max) / 2;
    ( t < xi[middle] ) ? (max = middle) : (min = middle);
  }
  // now xi[min] <= t < xi[min + 1]
  int32_t aux, count;
  const int32_t slack = 6;
  const int32_t group = 20;
  int64_t exp = n - 1;
  double prod = 1.0;
  // multiplying the factors to the left of t
  int32_t nLeft = min + 1;
  int32_t r = nLeft % group;
  if( r < slack ) {
    r = (nLeft > slack) ? (r + group) : nLeft;
  }

  for(int32_t j = 0; j < r; ++j) {
    prod = frexp(prod * (t - xi[min--]), &aux);
    exp  = exp + aux;
  }
  count = (min + 1) / group;
  for(int32_t i = 0; i < count; ++i) {
    for(int32_t j = 0; j < group; ++j) prod *= (t - xi[min--]);
    prod = frexp(prod, &aux);
    exp = exp + aux;
  }
  // multiplying the factors to the right of t
  int32_t nRight = n + 1 - max;
  r = nRight % group;
  if( r < slack ) {
    r = (nRight > slack) ? (r + group) : nRight;
  }

  for(int32_t j = 0; j < r; ++j) {
    prod = frexp(prod * (t - xi[max++]), &aux);
    exp  = exp + aux;
  }
  count = (n + 1 - max) / group;
  for(int32_t i = 0; i < count; ++i) {
    for(int32_t j = 0; j < group; ++j )  prod *= (t - xi[max++]);
    prod = frexp(prod, &aux);
    exp = exp + aux;
  }
  return scalbln(prod, exp);
}
\end{verbatim}
} 


\begin{thebibliography}{9}
\bibitem{BERRUT_POLE} J. --P. Berrut, 1988, Rational functions for guaranteed and experimentally well conditioned global interpolation, Comput. Math. Appl. 15 pp 1-16.
\bibitem{BERRUT} J. --P. Berrut and L. N. Trefethen, 2004, Barycentric Lagrange Interpolation, SIAM Review, 46, No. 3, pp. 501-517.
\bibitem{BOS} L. Bos, S. De Marchi and K. Hormann, 2011, On the Lebesgue constant of Berrut's rational interpolant at equidistant nodes, J. Comput. Appl. Math., 236(4), 504--510.
\bibitem{BOOST_SITE} The boost library website: http://www.boost.org/
\bibitem{BOOST_DEMING} R. Demming and D. Duffy, 2012, Introduction to the Boost C++ Libraries - Volume 2 - Advanced Libraries. Datasim. ISBN 978-94-91028-02-1.
\bibitem{CORMEN} T. H. Cormen et. al., 1990, Introduction to Algorithms (1st ed.). MIT Press and McGraw-Hill.
\bibitem{FLOATER} M. Floater and K. Hormann, 2007, Barycentric rational interpolation with no poles and high rates of approximation, Numer. Math. 107, 315-331.
\bibitem{FELLER} W. Feller, {\it An Introduction to Probability Theory and Its Applications}, Vol. 2, 3rd ed. New York: Wiley, 1971.
\bibitem{HENRICI} P. Henrici, 1964, {\it Elements of Numerical Analysis}, John Wiley and Sons, New York.
\bibitem{HIGHAM_SUM} N. J. Higham, 1993, The accuracy of floating point summation,
SIAM J. Scient. Comput. 14 (4) pp. 783--799.
\bibitem{HIGHAM} N.  J. Higham, 2004, The numerical stability of barycentric Lagrange
 interpolation, IMA J. Numer. Anal. 24, 547--556.
\bibitem{HIGHAM_BOOK} N.  J. Higham, 2002, {\it Accuracy and Stability of Numerical Algorithms}, 2nd ed., SIAM, Philadelphia.
\bibitem{KAHAN} W. Kahan, 1965, Further remarks on reducing truncation errors, Communications of the ACM 8 (1): 40.
\bibitem{RUMP} S. M. Rump, 2009, Ultimately Fast Accurate Summation, SIAM J. Scient. Comput. 31:5, pp. 3466--3502.
\bibitem{PRIESTA} D. Priest, 1991, Algorithms for arbitrary precision floating point arithmetic, in Proceedings of the
10th Symposium on Computer Arithmetic, Grenoble, France, P. Kornerup and D. Matula,
eds., IEEE Computer Society Press, Piscataway, NJ, pp. 132--145.
\bibitem{PRIESTB} D. Priest, 1992, On Properties of Floating Point Arithmetics: Numerical Stability and the Cost of
Accurate Computations, Ph.D. thesis, Mathematics Department, University of California at Berkeley, CA,
\bibitem{MPFR} L. Fousse, G. Hanrot, V. Lef\`{e}vre, P. P\'{e}lissier and P. Zimmermann, (2007),
MPFR: A Multiple-Precision Binary Floating-Point Library with Correct Rounding, ACM Transactions on Mathematical Software.
\bibitem{SALZER} H. Salzer, 1972, Lagrangian interpolation at the Chebyshev points $x_{n,\nu} \equiv \wfc{\cos}{\nu \pi/n},\nu = 0(1)n$;
some unnoted advantages, Comput. J., 15(2), pp. 156--159.
\bibitem{TAYLOR} W. J. Taylor, 1945, Method of Lagrangian curvilinear interpolation, J. Res. Nat. Bur. Standards, v. 35. pp 151-155.
\bibitem{ATAP} L. N. Trefethen, 2013, {\it Approximation Theory and Approximation Practice}, SIAM, Philadelphia.
\bibitem{WEBBX} M. Webb., L. N. Trefethen and P. Gonnet, 2012, Stability of barycentric interpolation formulas for extrapolation, SIAM J. Scient. Comput. 34 pp. 3009-3015.
\bibitem{WERNER} W. Werner, 1984, Polynomial interpolation: Lagrange versus Newton, Math. Comp. 43, pp 205-207.
\end{thebibliography}
\end{document}